\documentclass[review,onefignum,onetabnum]{siamonline220329}

\usepackage{lipsum}
\usepackage{amsfonts}
\usepackage{float}
\usepackage{graphicx}
\usepackage{caption}
\usepackage{epstopdf}
\usepackage{algorithmic}
\usepackage{amssymb}
\usepackage{amsmath}
\usepackage{tikz}
\usepackage{amsopn}

\nolinenumbers

\ifpdf
  \DeclareGraphicsExtensions{.eps,.pdf,.png,.jpg}
\else
  \DeclareGraphicsExtensions{.eps}
\fi

\usetikzlibrary{calc,shapes.geometric}
\headers{Convexification Method for Moving Targets}{M. V. Klibanov et al.}

\begin{document}


\title{Carleman Numerical Method for Imaging of Moving Targets 
	\thanks{Submitted to the editors DATE. }
}

\author{Michael V. Klibanov\thanks{%
Department of Mathematics and Statistics, University of North Carolina at
Charlotte, Charlotte, NC 28223, USA, (mklibanv@charlotte.edu)} 
\and Jingzhi
Li\thanks{%
Department of Mathematics, Southern University of Science and Technology,
Shenzhen 518055, P.~R.~China (li.jz@sustech.edu.cn)} \and Vladimir G. Romanov%
\thanks{%
Sobolev Institute of Mathematics, Novosibirsk 630090, Russian Federation
(romanov@math.nsc.ru)} \and Zhipeng Yang\thanks{%
School of Mathematics and Statistics, Lanzhou University, Lanzhou 730000, P.
R. China } }
\maketitle

\begin{abstract}
The problem of imaging of a moving target is formulated as a Coefficient
Inverse Problem for a hyperbolic equation with its coefficient depending on
all three spatial variables and time. As the initial condition, the point
source running along a straight line is used. Lateral Cauchy data are known
for each position of the point source. A truncated Fourier series with
respect to a special orthonormal basis is used. First, Lipschitz stability
estimate is obtained. Next, a globally convergent numerical method, the
so-called convexification method, is developed and its convergence analysis
is carried out. The convexification method is based on a Carleman estimate.
Results of numerical experiments are presented.
\end{abstract}



\section{ Introduction}

\label{sec:1}

Below $\mathbf{x=}\left( x,y,z\right) \in \mathbb{R}^{3}$ are points in $%
\mathbb{R}^{3}$ and $t\in \left( 0,T\right) $ is time, where $T>0$ is a
number. This is the first publication, in which the problem of imaging of
moving targets is posed as a full nonlinear Coefficient Inverse Problem
(CIP) of finding a generic $\left( \mathbf{x},t\right) -$dependent
coefficient of a hyperbolic equation. As the input data, we use Dirichlet
and Neumann boundary conditions, which are generated by the point source
running along an interval of a straight line. These input data are formally
determined ones, i.e. the number $m$ of free variables in the data equals
the number $n$ of free variables in the unknown coefficient, $m=n=4$.

CIPs of the determination of the $\left( \mathbf{x},t\right) -$dependent
coefficients using formally determined input data were not considered in the
past. Indeed, the majority of CIPs assume only the $\mathbf{x}-$dependence,
see, e.g. \cite{BukhKlib}, \cite{Isakov}-\cite{Kepid}. We also refer to,
e.g. \cite{Ali1,Ali2} and references cited therein for some CIPs with $%
\left( \mathbf{x},t\right) -$dependent unknown coefficients. The input data
in these publications are generated by the Dirichlet-to-Neumann maps with $%
m>n$. Uniqueness theorems were proven in these works, although numerical
methods were not developed.

A special case of interest for the moving target is when the unknown
coefficient $a\left( \mathbf{x},t\right) $ has the form $a\left( \mathbf{x}%
,t\right) =b\left( \mathbf{x}-\mathbf{v}t\right) ,$ where the function $b$
is unknown and $\mathbf{v}$ is the unknown velocity vector. However, rather
than using this assumption here, we find the coefficient $a\left( \mathbf{x}%
,t\right) $ in \textquotedblleft its entirety". In other words, we impose
almost minimal conditions on $a\left( \mathbf{x},t\right) .$ For each moment
of time $t$ we find location of the moving target, its shape and the spatial
distribution of this coefficient $a\left( \mathbf{x},t\right) ,$ both inside
and outside of that target.

\textbf{Definition.} \emph{We call a numerical method for a problem globally
convergent if there a theorem is proven, which claims that this method
provides points in a sufficiently small neighborhood of the true solution of
this problem without an advanced knowledge of that neighborhood.}

Our numerical method is the so-called convexification method, and it is
based on a Carleman estimate for a hyperbolic operator. The key advantage of
the convexification is its rigorously guaranteed global convergence, which
we establish for our CIP in our convergence analysis. Next, we provide
results of numerical experiments.

The problem of imaging of moving targets has an obvious applied interest.
One of methodologies used for this purpose is SAR imaging, which, however,
is based on the linearization via Born approximation, see, e.g. \cite%
{Borcea2013,Martorella}. Another interesting development in this direction
is imaging of moving targets in random media, see, e.g. \cite{Borcea2019,
Borcea20191} and references cited therein.

First, we reduce our CIP to a boundary value problem (BVP), which does not
contain the unknown coefficient. We prove Lipschitz stability estimate and
uniqueness for this problem. Next, we work on our main goal: the development
of the globally convergent convexification numerical method for the
reconstruction of the $\left( \mathbf{x},t\right) -$dependent coefficient.
To do this, we actually develop that method for that BVP.

Since our main goal is to develop a numerical method for our problem, then
we use here a simplification, assuming that the logarithm of the solution of
the forward problem admits an expansion in a truncated Fourier-like series
with respect to a special orthonormal basis, whose functions depend on the
position of the point source. This means that we come up with a new
mathematical model for our CIP. To our best knowledge, the idea of the
logarithmic transformation in CIPs was first published in \cite[second line
below formula (1.14)]{Klib95}. As to that basis, it was originally
introduced in \cite{Klib2017}. It turns out that the basis of \cite{Klib2017}
is a quite effective computational tool for a number of CIPs, see, e.g. \cite%
{KL,Trung} \ for some examples.

We now comment on the numerical method of this paper. This is a new version
of the above mentioned convexification method. Convexification was first
proposed in \cite{Klib95,Klib97} with the goal to avoid the well known
phenomenon of multiple local minima and ravines of conventional least
squares cost functionals for CIPs. These functionals are non convex. The
convexification method works for many CIPs. While publications \cite%
{Klib95,Klib97} are purely theoretical ones, more recent publications
contain both the theory and numerical studies, see, e.g. \cite%
{Bak,KL,Klibgrad,KMFG,Kepid} for some samples.

The convexification method is a numerical development of the idea of the
purely theoretical work \cite{BukhKlib}, in which the tool of Carleman
estimates was introduced in the field of Inverse Problems, also, see, e.g. 
\cite{Isakov,Klib92,Ksurvey,KL} for some follow up publications. For a given
CIP, the convexification constructs a weighted Tikhonov-like functional, in
which the weight is the Carleman Weight Function (CWF). This function is
used as the weight in the Carleman estimate for the corresponding PDE
operator. First, it is proven then that this functional is strongly convex
on an appropriate convex bounded set with an arbitrary diameter $d>0.$ Next,
the global convergence to the true solution of that CIP of the gradient
descent method of the minimization of this functional is proven and explicit
convergence rate is provided. Since no restrictions are imposed on the
number $d$, then this is global convergence by the above Definition.

An important part of the convexification method for any CIP\ is a
transformation procedure, which transforms that CIP in an ill-posed boundary
value problem (BVP)\ for either a single nonlinear integral differential PDE
or a system of coupled nonlinear PDEs. The unknown coefficient is not
present in both of these.\ In our case this is a system of coupled nonlinear
hyperbolic PDEs supplied by both Dirichlet and Neumann boundary conditions
on the lateral boundary of the time cylinder. These are the so-called
\textquotedblleft lateral Cauchy data". Initial conditions are not given in
this case. Thus, prior the construction of the above mentioned weighted
Tikhonov-like functional for this BVP, we prove Lipschitz stability and
uniqueness for it. The method of this proof is also based on a Carleman
estimate. This method was first proposed in \cite{KM} and was used since
then in a number of other publications, see, e.g. \cite{Ksurvey} and \cite[%
sections 2.7 and 3.6]{KL}.

\textbf{Remark 1.1.} Even though the convexification method works for many
CIPs, each new CIP requires its own convergence analysis due to its own
peculiarities, and this is what is done in the analytical part of the
current paper.

In section 2 we formulate both forward and inverse problems. In section 3 we
formulate a transformation procedure, which transforms the original CIP in a
boundary value problem for a system of coupled quasilinear hyperbolic
equations with the lateral Cauchy data. In section 4 we establish Lipschitz
stability and uniqueness for the above mentioned BVP. In section 5 we
present a version of the convexification method for that boundary value
problem and carry out convergence analysis. In section 6 we present results
of our numerical experiments. Finally, we prove one theorem in Appendix,
which is section 7. All functions below are real valued ones.

\section{Statements of Forward and Inverse Problems}

\label{sec:2}

Let $R>0$ be a number. Consider the ball $\Omega $ 
\begin{equation}
\Omega =\left\{ \mathbf{x=}\left( x,y,z\right) :\left\vert \mathbf{x}%
\right\vert <R\right\} .  \label{2.01}
\end{equation}%
For any number $T>0$ denote $D_{T}^{4}=\mathbb{R}^{3}\times \left(
0,T\right) .$ Set the positions of the point source 
\begin{equation}
\mathbf{x}_{0}=\left( s,0,-2R\right) ,s\in \left( -R,R\right) ,  \label{2.02}
\end{equation}%
where $s$ is a parameter. Hence, the point source runs along the interval $%
s\in \left( -R,R\right) $ of the line $L=\left\{ \mathbf{x}:x\in \mathbb{R}%
,y=0,z=-2R\right\} .$ Using (\ref{2.02}), denote the interval of our
interest on this line as 
\begin{equation}
L_{0}=\left\{ \mathbf{x}_{0}=\left( s,0,-2R\right) ,s\in \left( -R,R\right)
\right\} .  \label{2.03}
\end{equation}%
We take $-2R$ in (\ref{2.03}) rather than $-R$ to ensure that $\overline{%
\Omega }\cap L_{0}=\varnothing .$

Choose a number $T^{-}$ such that 
\begin{equation}
T>T^{-}>R\left( \sqrt{5}+1\right) .\text{ }  \label{2.04}
\end{equation}%
Denote 
\begin{equation}
T_{0}=\frac{T+T^{-}}{2},\text{ }A=\frac{T-T^{-}}{6}.  \label{1.40}
\end{equation}%
By (\ref{1.40}) 
\begin{equation}
T=T_{0}+3A,\text{ }T^{-}=T_{0}-3A.  \label{1.4}
\end{equation}%
Denote%
\begin{equation}
Q\left( R,T^{-},T\right) =\left\{ \left( \mathbf{x},t\right) :\mathbf{x}\in
\Omega ,t\in \left( T^{-},T\right) \right\} =\Omega \times \left(
T^{-},T\right) ,  \label{1.5}
\end{equation}%
\begin{equation}
S\left( R,T^{-},T\right) =\partial \Omega \times \left( T^{-},T\right)
=\left\{ \left( \mathbf{x},t\right) :\left\vert \mathbf{x}\right\vert
=R,t\in \left( T^{-},T\right) \right\} .  \label{1.50}
\end{equation}%
Hence, by (\ref{2.01}), (\ref{2.03}), (\ref{2.04}), (\ref{1.40})-(\ref{1.50})%
\begin{equation}
\left. 
\begin{array}{c}
\left\vert \mathbf{x-x}_{0}\right\vert <d=R\left( \sqrt{5}+1\right) ,\forall 
\mathbf{x}\in \Omega ,\text{ }\forall \mathbf{x}_{0}\in L_{0}, \\ 
t-\left\vert \mathbf{x-x}_{0}\right\vert >0,\text{ }\forall \left( \mathbf{x}%
,t\right) \in Q\left( R,T^{-},T\right) ,\text{ }\forall \mathbf{x}_{0}\in
L_{0}.%
\end{array}%
\right.  \label{1.7}
\end{equation}%
It follows from (\ref{1.5})-(\ref{1.7}) that 
\begin{equation}
Q\left( R,T^{-},T\right) ,\text{ }S\left( R,T^{-},T\right) \subset \left\{
t>\left\vert \mathbf{x-x}_{0}\right\vert \right\} ,\text{ }\forall \mathbf{x}%
_{0}\in L_{0}.  \label{1.8}
\end{equation}%
In other words, both sets $Q\left( R,T^{-},T\right) $ and $S\left(
R,T^{-},T\right) $ are located inside of the cone

$\left\{ t>\left\vert \mathbf{x-x}_{0}\right\vert \right\} $ for any
location of the source $\mathbf{x}_{0}\in L_{0}.$ Note that this is the
characteristic cone for the forward problem, which we formulate below.

Let the function $a\left( \mathbf{x},t\right) $ be the one characterizing
the moving target. We impose the following conditions on this function: 
\begin{equation}
a\left( \mathbf{x},t\right) \in C^{2}\left( \overline{D_{T}^{4}}\right) ,
\label{1.9}
\end{equation}%
\begin{equation}
a\left( \mathbf{x},t\right) \geq 0,  \label{1.10}
\end{equation}%
\begin{equation}
a\left( \mathbf{x},t\right) =0\text{ for }D_{T}^{4}\diagdown Q\left(
R,T^{-},T\right) .  \label{1.11}
\end{equation}

\textbf{The Forward Problem.} Suppose that the function $a\left( \mathbf{x}%
,t\right) $ is known and satisfies conditions (\ref{1.9})-(\ref{1.11}). For
each $\mathbf{x}_{0}\in L_{0},$ find the solution $u\left( \mathbf{x}%
,t\right) $ of the following Cauchy problem%
\begin{equation}
u_{tt}=\Delta u+a\left( \mathbf{x},t\right) u+\delta \left( \mathbf{x-x}%
_{0}\right) ,\left( \mathbf{x},t\right) \in D_{T}^{4},  \label{1.12}
\end{equation}%
\begin{equation}
u\left( \mathbf{x},0\right) =0,u_{t}\left( \mathbf{x},0\right) =0,
\label{1.13}
\end{equation}%
where $\delta \left( \mathbf{x-x}_{0}\right) $ is the delta function with
the position $\left\{ \mathbf{x}=\mathbf{x}_{0}\right\} $ of the point
source.

\textbf{Remark 2.1.} In principle, it would be better to replace equation (%
\ref{1.12}) with 
\begin{equation*}
c\left( \mathbf{x},t\right) u_{tt}=\Delta u+a\left( \mathbf{x},t\right)
u+\delta \left( \mathbf{x-x}_{0}\right) ,\left( \mathbf{x},t\right) \in
D_{T}^{4}
\end{equation*}%
with a positive function $c\left( \mathbf{x},t\right) $ and to assume in the
inverse problem that either of functions $c\left( \mathbf{x},t\right) $ or $%
a\left( \mathbf{x},t\right) $ is unknown whereas the other one is known.
However, due to the presence of the variable coefficient $c\left( \mathbf{x}%
,t\right) $ in the principal part of the hyperbolic operator, this case is
substantially more challenging than the one we consider here. Hence, on this
first stage of our research in this direction, it is reasonable to restrict
our attention to the case of equation (\ref{1.12}).

For each $\mathbf{x}_{0}\in L_{0}$ denote \ $\mathbf{x}_{0}\in L_{0}$%
\begin{equation}
K_{\mathbf{x}_{0},T}^{+}=\{(\mathbf{x},t)|\,|\mathbf{x}-\mathbf{x}_{0}|\leq
t\leq T\},  \label{200}
\end{equation}
Note that 
\begin{equation}
Q\left( R,T^{-},T\right) \subset \left( K_{\mathbf{x}_{0},T}^{+}\cap \left( 
\overline{\Omega }\times \left[ 0,T\right] \right) \right) ,\text{ }\forall 
\mathbf{x}_{0}\in L_{0}.  \label{201}
\end{equation}

Theorem 2.1 is proven in Appendix.

\textbf{Theorem 2.1.} \emph{Let conditions (\ref{1.9})-(\ref{1.11}) be
satisfied. Then the solution of problem (\ref{1.12}), (\ref{1.13}) has the
following form: }%
\begin{equation}
u(\mathbf{x},\mathbf{x}_{0},t)=H(t-|\mathbf{x}-\mathbf{x}_{0}|)\left[ \frac{1%
}{4\pi |\mathbf{x}-\mathbf{x}_{0}|}+\overline{u}(\mathbf{x},\mathbf{x}_{0},t)%
\right] ,\text{ }\forall \mathbf{x}_{0}\in L_{0},  \label{101}
\end{equation}%
\emph{where }$H(t)$\emph{\ is the Heaviside function} 
\begin{equation*}
H\left( t\right) =\left\{ 
\begin{array}{c}
1,t\geq 0, \\ 
0,t<0,%
\end{array}%
\right.
\end{equation*}%
\emph{and the function} $\overline{u}\in C^{2}\left( K_{\mathbf{x}%
_{0},T}^{+}\cap (\overline{\Omega }\times \lbrack 0,T])\right) $\emph{, }$%
\forall \mathbf{x}_{0}\in L_{0}.$\emph{\ Moreover,} 
\begin{equation}
u\left( \mathbf{x},\mathbf{x}_{0},t\right) \geq \frac{1}{4\pi |\mathbf{x}-%
\mathbf{x}_{0}|}\>\text{ for}\>\text{ }(\mathbf{x},t)\in K_{\mathbf{x}%
_{0},T}^{+}.  \label{102}
\end{equation}

It follows from (\ref{1.12}) and (\ref{1.13}) that 
\begin{equation}
u\left( \mathbf{x},\mathbf{x}_{0},t\right) =0\text{ for }t\in \left(
0,\left\vert \mathbf{x}-\mathbf{x}_{0}\right\vert \right) .  \label{1.16}
\end{equation}%
This explains why do we need the second line of (\ref{1.7}) in combination
with (\ref{1.4}), (\ref{1.5}) and (\ref{1.11}). Let $\nu \left( \mathbf{x}%
\right) ,\mathbf{x}\in \partial \Omega $ be the outward normal vector on $%
\partial \Omega $ and $\partial _{\nu }$ be the normal derivative.

\textbf{Coefficient Inverse Problem. }\emph{Let }$S\left( R,T^{-},T\right) $%
\emph{\ be the cylindrical hypersurface defined in (\ref{1.50}). Assume that
the following two functions }$g_{0}\left( \mathbf{x},\mathbf{x}_{0},t\right)
,g_{1}\left( \mathbf{x},\mathbf{x}_{0},t\right) $\emph{\ are known:}%
\begin{equation}
\left. 
\begin{array}{c}
u\left( \mathbf{x},\mathbf{x}_{0},t\right) \mid _{S\left( R,T^{-},T\right)
}=g_{0}\left( \mathbf{x},\mathbf{x}_{0},t\right) ,\text{ } \\ 
\partial _{\nu }u\left( \mathbf{x},\mathbf{x}_{0},t\right) \mid _{S\left(
R,T^{-},T\right) }=g_{1}\left( \mathbf{x},\mathbf{x}_{0},t\right) , \\ 
\forall \mathbf{x}_{0}\in L_{0}.%
\end{array}%
\right.  \label{1.190}
\end{equation}%
\emph{\ Assume that the coefficient }$a\left( \mathbf{x},t\right) $\emph{\
satisfies conditions (\ref{1.9})-(\ref{1.11}) and it is unknown. Find the
function \ }$a\left( \mathbf{x},t\right) $\emph{\ for} $\left( \mathbf{x}%
,t\right) \in Q\left( R,T^{-},T\right) .$

\textbf{Physical Meaning.} In acoustics 
\begin{equation*}
a\left( \mathbf{x}\right) =\frac{\Delta \left( \sqrt{\rho \left( \mathbf{x}%
\right) }\right) }{\sqrt{\rho \left( \mathbf{x}\right) }},
\end{equation*}%
where $\rho \left( \mathbf{x}\right) >0$ is the density of the medium \cite%
{Rom1}. The physical meaning of problem (\ref{1.12}), (\ref{1.13}) is that
the point source located at $\left\{ \mathbf{x}_{0}\right\} $ radiates the
wave field $u\left( \mathbf{x},\mathbf{x}_{0},t\right) $ for all times $t>0.$
In the inverse problem, one measures the scattered wave field and its normal
derivative on the surface of the ball $\Omega $ for all times $t\in \left(
T^{-},T\right) .$ Using these measurements, one wants to recover the
coefficient $a\left( \mathbf{x},t\right) .$ This coefficient represents the
unknown moving target. The fact that the unknown target is moving, is
reflected by the time dependence of $a\left( \mathbf{x},t\right) .$

\section{Transformation}

\label{sec:3}

\subsection{An equation without $a\left( \mathbf{x},t\right) $}

\label{sec:3.1}

It follows from (\ref{102}) that we can consider in

$Q\left( R,T^{-},T\right) $ the function $v\left( \mathbf{x},\mathbf{x}%
_{0},t\right) ,$%
\begin{equation}
v\left( \mathbf{x},\mathbf{x}_{0},t\right) =\ln \left[ u\left( \mathbf{x},%
\mathbf{x}_{0},t\right) \right] \in C^{2}\left( \overline{Q\left(
R,T^{-},T\right) }\right) ,\forall \mathbf{x}_{0}\in L_{0}.  \label{2.4}
\end{equation}%
Since by (\ref{2.01})-(\ref{2.03}) and (\ref{1.5}) $L_{0}\cap \overline{%
Q\left( R,T^{-},T\right) }=\varnothing ,$ then 
\begin{equation}
\delta \left( \mathbf{x-x}_{0}\right) =0\text{, }\forall \mathbf{x}\in 
\overline{\Omega },\text{ }\forall \mathbf{x}_{0}\in L_{0}.  \label{2.5}
\end{equation}%
Hence, substituting (\ref{2.4}) in equation (\ref{1.12}) without the term in
(\ref{2.5}), we obtain%
\begin{equation}
v_{tt}-\Delta v+v_{t}^{2}-\left( \nabla v\right) ^{2}=a\left( \mathbf{x}%
,t\right) ,\text{ }\left( \mathbf{x},t\right) \in \text{ }Q\left(
R,T^{-},T\right) ,\text{ }\mathbf{x}_{0}\in L_{0}.  \label{2.6}
\end{equation}%
In addition, by (\ref{2.03}) and (\ref{1.190})%
\begin{equation}
\left. 
\begin{array}{c}
v\left( \mathbf{x},\mathbf{x}_{0},t\right) \mid _{\left( \mathbf{x},t\right)
\in S\left( R,T^{-},T\right) }=p_{0}\left( \mathbf{x},\mathbf{x}%
_{0},t\right) , \\ 
\partial _{\nu }v\left( \mathbf{x},\mathbf{x}_{0},t\right) \mid _{\left( 
\mathbf{x},t\right) \in S\left( R,T^{-},T\right) }=p_{1}\left( \mathbf{x},%
\mathbf{x}_{0},t\right) , \\ 
p_{0}\left( \mathbf{x},\mathbf{x}_{0},t\right) =\ln \left[ g_{0}\left( 
\mathbf{x},\mathbf{x}_{0},t\right) \right] ,\text{ }p_{1}\left( \mathbf{x},%
\mathbf{x}_{0},t\right) =\left( g_{1}/g_{0}\right) \left( \mathbf{x},\mathbf{%
\ x}_{0},t\right) , \\ 
\mathbf{x}_{0}=\left( s,0,-2R\right) \in L_{0},\text{ }s\in \left(
-R,R\right) .%
\end{array}%
\right.  \label{2.7}
\end{equation}%
Differentiate (\ref{2.6}) and (\ref{2.7}) with respect to the parameter $%
s\in \left( -R,R\right) $ in (\ref{2.7}) and use 
\begin{equation*}
\frac{\partial a\left( \mathbf{x},t\right) }{\partial s}\equiv 0.
\end{equation*}%
We obtain%
\begin{equation}
\left. 
\begin{array}{c}
\partial _{t}^{2}v_{s}-\Delta v_{s}+2\partial _{t}\left( v_{s}\right)
\partial _{t}v-2\nabla v_{s}\cdot \nabla v=0\text{ in }Q\left(
R,T^{-},T\right) ,s\in \left( -R,R\right) , \\ 
v\left( \mathbf{x},\mathbf{x}_{0},t\right) \mid _{S\left( R,T^{-},T\right)
}=p_{0}\left( \mathbf{x},\mathbf{x}_{0},t\right) , \\ 
\partial _{\nu }v\left( \mathbf{x},\mathbf{x}_{0},t\right) \mid _{S\left(
R,T^{-},T\right) }=p_{1}\left( \mathbf{x},\mathbf{x}_{0},t\right) , \\ 
v_{s}\left( \mathbf{x},\mathbf{x}_{0},t\right) \mid _{S\left(
R,T^{-},T\right) }=\partial _{s}p_{0}\left( \mathbf{x},\mathbf{x}%
_{0},t\right) , \\ 
\partial _{\nu }v_{s}\left( \mathbf{x},\mathbf{x}_{0},t\right) \mid
_{S\left( R,T^{-},T\right) }=\partial _{s}p_{1}\left( \mathbf{x},\mathbf{x}%
_{0},t\right) , \\ 
\forall \mathbf{x}_{0}=\left( s,0,-2R\right) ,\forall s\in \left(
-R,R\right) .%
\end{array}%
\right.  \label{2.9}
\end{equation}

Below we focus on finding functions $v,v_{s}$ from conditions (\ref{2.9}).
Suppose that a solution of this problem is found. In the case of
computations, denote the latter function $v_{\text{comp}}\left( \mathbf{x}%
,s,t\right) .$ Then we substitute this function in equation (\ref{2.6}) and
set%
\begin{equation}
\left. 
\begin{array}{c}
a_{\text{comp}}\left( \mathbf{x},t\right) = \\ 
=\left( 1/2R\right) \int\limits_{-R}^{R}\left[ \partial _{t}^{2}v_{\text{comp%
}}-\Delta v_{\text{ comp}}+\left( \partial _{t}v_{\text{comp}}\right)
^{2}-\left( \nabla v_{\text{comp}}\right) ^{2}\right] \left( \mathbf{x}%
,s,t\right) ds, \\ 
\left( \mathbf{x},t\right) \in Q\left( R,T^{-},T\right) .%
\end{array}%
\right.  \label{2.10}
\end{equation}

\subsection{Orthonormal basis}

\label{sec:3.2}

It is hard to solve problem (\ref{2.9}) in its present form. Hence, keeping
in mind that the main goal of this paper is a version of the convexification
numerical method rather than a purely theoretical result, we simplify our
task via introduction of a new mathematical model, which uses a special
orthonormal basis. This basis was first introduced in \cite{Klib2017}, also
see \cite[section 6.2.3]{KL}.

Consider the set of functions $\left\{ s^{n}e^{s}\right\} _{n=0}^{\infty
},s\in \left( -R,R\right) .$ These functions are linearly independent and
form a complete set in $L_{2}\left( -R,R\right) .$

\textbf{Step 1.} Select an integer $N>1.$ Use the classical Gram-Schmidt
orthonormalization procedure to orthonormalize first $N$ elements of this
system. The number $N$\ needs to be chosen computationally.

\textbf{Step 2. }Therefore, we obtain functions 
\begin{equation}
\left\{ \psi _{n}\left( s\right) \right\} _{n=0}^{N-1},\text{ }\psi
_{n}\left( s\right) =P_{n}\left( s\right) e^{s},  \label{2.100}
\end{equation}%
where $P_{n}\left( s\right) $ is a polynomial of the degree $n.$

Obviously, the set $\left\{ \psi _{n}\left( s\right) \right\} _{n=0}^{\infty
}$ is an orthonormal basis in the space $L_{2}\left( -R,R\right) .$ Theorem
3.1 is the key result for this construction.

\textbf{Theorem} \textbf{3.1} (\cite{Klib2017}, \cite[section 6.2.3]{KL}). 
\emph{Let }$\left( ,\right) $\emph{\ be the scalar product in the space }$%
L_{2}\left( -R,R\right) .$\emph{\ Denote }%
\begin{equation}
a_{mk}=\left( \psi _{k}^{\prime },\psi _{m}\right) =\int\limits_{-R}^{R}\psi
_{k}^{\prime }\left( s\right) \psi _{m}\left( s\right) ds.  \label{3.1}
\end{equation}%
\emph{Then}%
\begin{equation}
a_{mk}=\left\{ 
\begin{array}{c}
1\text{ if }k=m, \\ 
0\text{ if }k<m.%
\end{array}%
\right.  \label{3.2}
\end{equation}%
\emph{Consider the matrix }$M_{N}=\left( a_{mk}\right) _{m,k=0}^{N-1}.$\emph{%
\ \ Then (\ref{3.2}) implies that }$\det M_{N}=1.$\emph{\ Therefore, there
exists the matrix }$M_{N}^{-1}.$

We note that the matrix $M_{N}$ is not invertible for both classical
orthonormal polynomials and functions forming the trigonometric basis. This
is because the function representing the identical constant is a part of
these bases.

\subsection{A boundary value problem for a system of nonlinear hyperbolic
equations}

\label{sec:3.3}

Consider an approximation of the function $v\left( \mathbf{x},s,t\right) $
in (\ref{2.6}) via the truncated Fourier-like series and replace in this
approximation the symbol \textquotedblleft $\approx "$ with the symbol
\textquotedblleft $=$",%
\begin{equation}
v\left( \mathbf{x},s,t\right) =\sum\limits_{n=0}^{N-1}v_{k}\left( \mathbf{x}%
,t\right) \psi _{k}\left( s\right) ,\text{ }\left( \mathbf{x},t\right) \in 
\text{ }Q\left( R,T^{-},T\right) ,s\in \left( -R,R\right) ,  \label{4.1}
\end{equation}%
where Fourier coefficients $\left\{ v_{n}\left( \mathbf{x},t\right) \right\}
_{n=0}^{N-1}$ are unknown and are the target of our search. Then%
\begin{equation}
\partial _{s}v\left( \mathbf{x},s,t\right)
=\sum\limits_{n=0}^{N-1}v_{k}\left( \mathbf{x},t\right) \psi _{k}^{\prime
}\left( s\right) ,\text{ }\left( \mathbf{x},t\right) \in Q\left(
R,T^{-},T\right) ,s\in \left( -R,R\right) .  \label{4.2}
\end{equation}%
Substitute (\ref{4.1}) and (\ref{4.2}) in the first line of (\ref{2.9}). We
obtain%
\begin{equation}
\left. 
\begin{array}{c}
\sum\limits_{n=0}^{N-1}\left( \partial _{t}^{2}v_{n}\mathbf{-}\Delta
v_{n}\right) \left( \mathbf{x},t\right) \psi _{k}^{\prime }\left( s\right)
+2\sum\limits_{n,k=0}^{N-1}\left( \partial _{t}v_{n}\partial
_{t}v_{k}\right) \left( \mathbf{x},t\right) \psi _{k}^{\prime }\left(
s\right) \psi _{n}\left( s\right) - \\ 
-2\sum\limits_{n,k=0}^{N-1}\nabla v_{n}\cdot \nabla v_{k}\psi _{k}^{\prime
}\left( s\right) \psi _{n}\left( s\right) =0, \\ 
\left( \mathbf{x},t\right) \in Q\left( R,T^{-},T\right) ,s\in \left(
-R,R\right) .%
\end{array}%
\right.  \label{4.3}
\end{equation}%
Consider the $N-$dimensional vector function $V\left( \mathbf{x},t\right) $
of coefficients $v_{n}\left( \mathbf{x},t\right) ,$%
\begin{equation}
V\left( \mathbf{x},t\right) =\left( v_{0},...,v_{N-1}\right) ^{T}\left( 
\mathbf{x},t\right) ,\text{ }\left( \mathbf{x},t\right) \in Q\left(
R,T^{-},T\right) .  \label{4.4}
\end{equation}%
Sequentially multiply equation (\ref{4.3}) by functions $\psi _{m}\left(
s\right) ,m=0,...,N-1$ and integrate each such product with respect to $s\in
\left( -R,R\right) .$ Then, using (\ref{3.1}), (\ref{3.2}) and (\ref{4.4}),
we obtain the following system of nonlinear hyperbolic equations with
respect to the vector function $V\left( \mathbf{x},t\right) :$%
\begin{equation}
M_{N}\left( V_{tt}-\Delta V\right) +F\left( \nabla V,V\right) =0,\text{ }%
\left( \mathbf{x},t\right) \in Q\left( R,T^{-},T\right) ,  \label{4.5}
\end{equation}%
where $\nabla V\left( \mathbf{x},t\right) =\left( \nabla v_{0},...,\nabla
v_{N-1}\right) ^{T}\left( \mathbf{x},t\right) $ and the $N-$D vector
function $F\left( \alpha \right) $ is 
\begin{equation}
F\left( \alpha \right) =\left( F^{\left( 0\right) },...,F^{\left( N-1\right)
}\right) ^{T}\left( \alpha \right) \in C^{1}\left( \mathbb{R}^{4N}\right)
,\alpha =\left( \alpha _{1},...,\alpha _{4N}\right) ^{T}\in \mathbb{R}^{4N}.
\label{4.06}
\end{equation}

The Dirichlet and Neumann boundary conditions for the vector function $%
V\left( \mathbf{x},t\right) $ at the hypersurface $S\left( R,T^{-},T\right) $
are:%
\begin{equation}
V\mid _{\left( \mathbf{x},t\right) \in S\left( R,T^{-},T\right)
}=q_{0}\left( \mathbf{x},t\right) ,\text{ }\partial _{\nu }V\mid _{\left( 
\mathbf{x},t\right) \in S\left( R,T^{-},T\right) }=q_{1}\left( \mathbf{x}
,t\right) .  \label{4.7}
\end{equation}

Our focus below is on the solution of the following Boundary Value Problem:

\textbf{Boundary Value Problem (BVP).} \emph{Find the vector function }$%
V\left( \mathbf{x},t\right) $\emph{\ for }

$\left( \mathbf{x},t\right) \in Q\left( R,T^{-},T\right) $ \emph{\ from
conditions (\ref{4.5})-(\ref{4.7}).}

If this problem is solved, then the function $v_{\text{comp}}\left( \mathbf{x%
},s,t\right) $ is found via (\ref{4.1}), (\ref{4.4}). Next, formula (\ref%
{2.10}) should be used to find the unknown coefficient $a\left( \mathbf{x}%
,t\right) .$

Let $M_{1}>0$ be a number. We seek such a solution of problem (\ref{4.5}), (%
\ref{4.06}) on the set of $N-$D vector functions $V\left( \mathbf{x}%
,t\right) \in C^{2}\left( \overline{Q\left( R,T^{-},T\right) }\right) $ in (%
\ref{4.4}) that 
\begin{equation}
\max_{\left( \mathbf{x},t\right) \in \overline{Q\left( R,T^{-},T\right) }%
}\left\vert \nabla V\left( \mathbf{x},t\right) \right\vert ,\max_{\left( 
\mathbf{x},t\right) \in \overline{Q\left( R,T^{-},T\right) }}\left\vert
V\left( \mathbf{x},t\right) \right\vert \leq M_{1}.  \label{4.8}
\end{equation}%
It follows from (\ref{4.06}) that there exists a constant $M_{2}=M_{2}\left(
F,M_{1}\right) >0$ such that%
\begin{equation}
\max_{\left\vert \alpha _{1}\right\vert ,...,\left\vert \alpha
_{4N}\right\vert \leq M_{1}}\left\vert \nabla _{\alpha }F\left( \alpha
\right) \right\vert \leq M_{2}.  \label{4.9}
\end{equation}

We now explain how to obtain vector functions $q_{0}\left( \mathbf{x}%
,t\right) $ and $q_{1}\left( \mathbf{x},t\right) $ in (\ref{4.7}). Using (%
\ref{2.7}), we obtain%
\begin{equation*}
\left. 
\begin{array}{c}
v_{k}\mid _{\left( \mathbf{x},t\right) \in S\left( R,T^{\pm }\right)
}=\int\limits_{-R}^{R}p_{0}\left( \mathbf{x},s,t\right) \psi _{k}\left(
s\right) ds=q_{0,k}\left( \mathbf{x},t\right) ,\text{ }\left( \mathbf{x}%
,t\right) \in S\left( R,T^{-},T\right) , \\ 
\partial _{\nu }v_{k}\mid _{\left( \mathbf{x},t\right) \in S\left( R,T^{\pm
}\right) }=\int\limits_{-R}^{R}p_{1}\left( \mathbf{x},s,t\right) \psi
_{k}\left( s\right) ds=q_{1,k}\left( \mathbf{x},t\right) ,\left( \mathbf{x}%
,t\right) \in S\left( R,T^{-},T\right) , \\ 
k=0,...,N-1.%
\end{array}%
\right.
\end{equation*}%
Next, we set 
\begin{equation*}
\left. 
\begin{array}{c}
q_{0}\left( \mathbf{x},t\right) =\left( q_{0,1},...,q_{0,N-1}\right)
^{T}\left( \mathbf{x},t\right) ,q_{1}\left( \mathbf{x},t\right) =\left(
q_{1,1},...,q_{1,N-1}\right) ^{T}\left( \mathbf{x},t\right) , \\ 
\left( \mathbf{x},t\right) \in S\left( R,T^{-},T\right) .%
\end{array}%
\right.
\end{equation*}%
The transformation procedure is completed.

\section{Lipschitz Stability and Uniqueness for Problem (\protect\ref{4.5}%
)-( \protect\ref{4.7})}

\label{sec:4}

Problem (\ref{4.5})-(\ref{4.7}) is the problem for the system (\ref{4.5}) of
coupled quasilinear hyperbolic PDEs with the lateral Cauchy data (\ref{4.7}%
). To obtain the desired Lipschitz stability estimate, we use the Carleman
estimate for the wave operator $\partial _{t}^{2}-\Delta $ \cite[Theorem
2.5.1]{KL}. To take into account conditions (\ref{2.01})-(\ref{1.8}), we
specify values of parameters in that estimate.

\subsection{Carleman estimate}

\label{sec:4.1}

Let $\sigma $,$\eta >0$ be two positive numbers, which we will specify
later. Let the point 
\begin{equation}
\mathbf{p}=\left( -R-\sigma ,0,0\right) .  \label{5.001}
\end{equation}%
Then by (\ref{2.01}), (\ref{5.001}) and triangle inequality%
\begin{equation}
\left\vert \mathbf{x}-\mathbf{p}\right\vert \in \left[ \sigma ,2R+\sigma %
\right] ,\text{ }\forall \mathbf{x}\in \overline{\Omega }.  \label{5.1}
\end{equation}%
Let a sufficiently small number $h$ be such that 
\begin{equation}
h\in \left( 0,\frac{\sigma ^{2}}{5}\right) .  \label{5.01}
\end{equation}%
Recalling (\ref{1.40}), we choose the number $T$ so large that the number $%
\eta $ satisfies: 
\begin{equation}
\eta =\frac{\sigma ^{2}-h}{\left( 2A\right) ^{2}}\in \left( 0,1\right) .
\label{5.2}
\end{equation}%
Consider the function $\psi \left( \mathbf{x},t\right) ,$%
\begin{equation}
\psi \left( \mathbf{x},t\right) =\left\vert \mathbf{x}-\mathbf{p}\right\vert
^{2}-\eta \left( t-T_{0}\right) ^{2}.  \label{5.3}
\end{equation}%
Consider the domain $G_{h},$%
\begin{equation}
G_{h}=\left\{ \left( \mathbf{x},t\right) :\mathbf{x}\in \Omega ,\text{ }\psi
\left( \mathbf{x},t\right) =\left\vert \mathbf{x}-\mathbf{p}\right\vert
^{2}-\eta \left( t-T_{0}\right) ^{2}>h\right\} .  \label{5.4}
\end{equation}

We choose the number $\sigma >0$ such that 
\begin{equation}
\left( 2R+\sigma \right) ^{2}-\frac{9}{4}\sigma ^{2}<-\frac{1}{4}\sigma ^{2}.
\label{5.60}
\end{equation}%
Since (\ref{5.01}) implies that $\sigma ^{2}>5h$, then the choice (\ref{5.60}%
) implies 
\begin{equation}
\left( 2R+\sigma \right) ^{2}-\frac{9\sigma ^{2}}{4}<-\frac{5}{4}h,
\label{5.61}
\end{equation}%
It follows from (\ref{5.2}) that (\ref{5.61}) is equivalent with 
\begin{equation}
\left( 2R+\sigma \right) ^{2}-9A^{2}\eta =\left( 2R+\sigma \right)
^{2}-9A^{2}\frac{\sigma ^{2}-h}{4A^{2}}<h.  \label{5.62}
\end{equation}%
Inequality (\ref{5.60}) is equivalent with%
\begin{equation*}
\left( 2R+\left( 1-\sqrt{2}\right) \sigma \right) \left( 2R+\left( 1+\sqrt{2}%
\right) \sigma \right) <0.
\end{equation*}%
\textbf{\ }Hence, to ensure (\ref{5.60}), we choose the number $\sigma $
such that%
\begin{equation}
\sigma >\frac{2R}{\sqrt{2}-1}.  \label{5.6}
\end{equation}%
It follows from (\ref{1.4}), (\ref{5.4}) and (\ref{5.62}) that (\ref{5.6})
implies 
\begin{equation*}
\text{ }\psi \left( \mathbf{x},T\right) <h,\text{ }\psi \left( \mathbf{x}%
,T^{-}\right) <h.
\end{equation*}%
Hence, 
\begin{equation}
\left. 
\begin{array}{c}
G_{h}\subset Q\left( R,T^{-},T\right) , \\ 
\text{ }\overline{G}_{h}\cap \left\{ t=T^{-}\right\} =\overline{G}_{h}\cap
\left\{ t=T\right\} =\varnothing .%
\end{array}%
\right.  \label{5.8}
\end{equation}%
By (\ref{5.1})-(\ref{5.2}) and (\ref{5.4})%
\begin{equation}
Q\left( R,A\right) =\Omega \times \left( T_{0}-A,T_{0}+A\right) \subset
G_{4h}\subset G_{3h}\subset G_{2h}\subset G_{h}.  \label{5.80}
\end{equation}

The boundary $\partial G_{h}$ of the domain $G_{h}$ consists of two parts:%
\begin{equation}
\left. 
\begin{array}{c}
\partial G_{h}=\partial _{1}G_{h}\cup \partial _{2}G_{h}, \\ 
\partial _{1}G_{h}=\left\{ \left( \mathbf{x},t\right) :\psi \left( \mathbf{x}
,t\right) =h\right\} \cap \overline{Q\left( R,T^{-},T\right) }, \\ 
\partial _{2}G_{h}=S\left( R,T^{-},T\right) \cap \overline{G}_{h}.%
\end{array}
\right.  \label{5.81}
\end{equation}

Let $\lambda >1$ be a large parameter, which we will specify later. Using (%
\ref{5.3}), we define the Carleman Weight Function (CWF) as%
\begin{equation}
\varphi \left( \mathbf{x},t\right) =e^{\lambda \psi \left( \mathbf{x}%
,t\right) }.  \label{5.9}
\end{equation}%
Theorem 2.5.1 of \cite{KL} implies that the following pointwise Carleman
estimate is valid:

\textbf{Theorem 4.1.} \emph{Assume that conditions (\ref{5.001}), (\ref{5.01}
), (\ref{5.2}), (\ref{5.3}), (\ref{5.4}), (\ref{5.6}) and (\ref{5.9}) hold.
Then there exists a sufficiently large number }

$\lambda _{0}=\lambda _{0}\left( \sigma ,h,T_{0},T,R,G_{h}\right) \geq 1$%
\emph{\ and a number }$C=C\left( \sigma ,h,T_{0},T,R,G_{h}\right) >0,$\emph{%
\ both numbers depending only on listed parameters, such that the following
pointwise Carleman estimate for the wave operator }$\partial _{t}^{2}-\Delta 
$\emph{\ holds:} 
\begin{equation}
\left( u_{tt}-\Delta u\right) ^{2}\varphi ^{2}\geq C\lambda \left(
u_{t}^{2}+\left\vert \nabla u\right\vert ^{2}\right) \varphi ^{2}+C\lambda
^{3}u^{2}\varphi ^{2}+\text{div}U+W_{t}\text{ in }G_{h},  \label{5.10}
\end{equation}%
\begin{equation}
\forall u\in C^{2}\left( \overline{G}_{h}\right) ,\text{ }\forall \lambda
\geq \lambda _{0},  \label{5.11}
\end{equation}%
\emph{where the vector function }$\left( U,W\right) =\left(
U_{1},U_{2},U_{3},W\right) \left( \mathbf{x},t\right) $\emph{\ satisfies the
following estimate:}%
\begin{equation}
\left\vert \left( U,W\right) \right\vert \left( \mathbf{x},t\right) \leq
C\lambda ^{3}\left( u_{t}^{2}+\left\vert \nabla u\right\vert
^{2}+u^{2}\right) \left( \mathbf{x},t\right) \varphi ^{2}\left( \mathbf{x}%
,t\right) ,\text{ }\forall \left( \mathbf{x},t\right) \in \overline{G}_{h}.
\label{5.12}
\end{equation}

\subsection{Lipschitz stability and uniqueness of BVP (\protect\ref{4.5})-( 
\protect\ref{4.7})}

\label{sec:4.2}

Let $B$ be a Banach space of functions with its norm $\left\Vert \cdot
\right\Vert _{B}.$ Let $B_{k}=\underbrace{B\times B\times ...\times B},$ $k$
times. To simplify notations, we use below for the norm in $B_{k}$ the same
notation $\left\Vert \cdot \right\Vert _{B}$ as the one for the norm in $B$,
meaning that%
\begin{equation*}
\left\Vert f\right\Vert _{B_{k}}=\left( \sum\limits_{i=1}^{k}\left\Vert
f_{i}\right\Vert _{B}^{2}\right) ^{1/2},\text{ }\forall f=\left(
f_{1},...,f_{k}\right) ^{T}\in B_{k}.
\end{equation*}

\textbf{Theorem 4.2.} \emph{Assume that conditions (\ref{2.01})-(\ref{1.50}
), (\ref{4.06}) and (\ref{4.9}) hold. Suppose that there exist two }$N-$%
\emph{D vector functions }$V^{\left( 1\right) }\left( \mathbf{x},t\right)
,V^{\left( 2\right) }\left( \mathbf{x},t\right) \in C_{N}^{2}\left( 
\overline{Q\left( R,T^{-},T\right) }\right) $\emph{\ satisfying equation (%
\ref{4.5}), condition (\ref{4.8}) and the following analogs of Dirichlet and
Neumann boundary conditions (\ref{4.7}):}%
\begin{equation}
V^{\left( i\right) }\mid _{S\left( R,T^{-},T\right) }=q_{0}^{\left( i\right)
}\left( \mathbf{x},t\right) ,\text{ }\partial _{\nu }V^{\left( i\right)
}\mid _{S\left( R,T^{-},T\right) }=q_{1}^{\left( i\right) }\left( \mathbf{x}%
,t\right) ,\text{ }i=1,2.  \label{5.13}
\end{equation}%
\emph{Denote} 
\begin{equation}
\widetilde{V}=V^{\left( 1\right) }-V^{\left( 2\right) },\widetilde{q}%
_{0}=q_{0}^{\left( 1\right) }-q_{0}^{\left( 2\right) },\widetilde{q}%
_{1}=q_{1}^{\left( 1\right) }-q_{1}^{\left( 2\right) }.  \label{5.14}
\end{equation}%
\emph{Suppose that conditions \ref{5.001}), (\ref{5.01}), (\ref{5.2}), (\ref%
{5.3}), (\ref{5.4}), (\ref{5.6}) and (\ref{5.9}) are valid. Then there
exists a number }

$C_{1}=C_{1}\left( Q\left( R,T^{-},T\right) ,M_{1},M_{2},N,h\right) >0$\emph{%
\ \ depending only on listed parameters but independent on }$\widetilde{V}$%
\emph{\ such that the following Lipschitz stability estimate holds:}%
\begin{equation}
\left\Vert \widetilde{V}\right\Vert _{H_{N}^{1}\left( Q\left(
R,T^{-},T\right) \right) }\leq C_{1}\left( \left\Vert \widetilde{q}%
_{0}\right\Vert _{H_{N}^{1}\left( S\right) }+\left\Vert \widetilde{q}%
_{1}\right\Vert _{L_{2,N}\left( S\left( R,T^{\pm }\right) \right) }\right) .
\label{5.15}
\end{equation}%
\emph{In particular, if }$\widetilde{q}_{0}\equiv \widetilde{q}_{1}\equiv 0,$%
\emph{\ then }$\widetilde{V}\equiv 0,$ \emph{which implies uniqueness of BVP
(\ref{4.5})-(\ref{4.7}). }

\textbf{Proof}. Below $C>0$ and $C_{1}>0$ denote different numbers depending
only on parameters listed in formulations of Theorems 4.1 and 4.2
respectively. Using Theorem 3.1, multiply from the left both sides of
equation (\ref{4.5}) by the matrix $M_{N}^{-1}$. Then subtract the obtained
vector equation for $V^{\left( 2\right) }$ from the obtained vector equation
for $V^{\left( 1\right) }.$ Next, using (\ref{4.06}), (\ref{4.8}) and (\ref%
{4.9}), (\ref{5.13}), (\ref{5.14}) and the multidimensional analog of the
Taylor formula (see, e.g. \cite{V} for that formula), we obtain a vector
equation for $\widetilde{V},$%
\begin{equation}
L_{0}\left( \widetilde{V}\right) =\widetilde{V}_{tt}-\Delta \widetilde{V}%
=B_{1}\left( \mathbf{x},t\right) \cdot \nabla \widetilde{V}+B_{2}\left( 
\mathbf{x},t\right) \widetilde{V}\text{ in }Q\left( R,T^{-},T\right) ,
\label{5.16}
\end{equation}%
\begin{equation}
\widetilde{V}\mid _{S\left( R,T^{-},T\right) }=\widetilde{q}_{0},\text{ }%
\partial _{\nu }\widetilde{V}\mid _{S\left( R,T^{-},T\right) }=\widetilde{q}%
_{1},  \label{5.160}
\end{equation}%
where $B_{1}\left( \mathbf{x},t\right) $ is a $N\times 3N$ matrix, $%
B_{2}\left( \mathbf{x},t\right) $ is a $N\times N$ matrix with their
elements belonging to the space $C\left( \overline{Q\left( R,T^{-},T\right) }%
\right) ,$ and 
\begin{equation}
\left\vert B_{1}\left( \mathbf{x},t\right) \right\vert ,\left\vert
B_{2}\left( \mathbf{x},t\right) \right\vert \leq C_{1},\text{ }\left( 
\mathbf{x},t\right) \in \overline{Q\left( R,T^{-},T\right) }.  \label{5.17}
\end{equation}%
The latter inequality means that $C\left( \overline{Q\left( R,T^{-},T\right) 
}\right) -$norms of all elements of both these matrices are bounded by $%
C_{1}.$

Based on (\ref{5.8}) and (\ref{5.80}), choose a cut-off function $\chi
\left( \mathbf{x},t\right) ,$%
\begin{equation}
\left. 
\begin{array}{c}
\chi \in C^{2}\left( \overline{Q\left( R,T^{-},T\right) }\right) ,\text{ }
\\ 
\chi \left( \mathbf{x},t\right) =\left\{ 
\begin{array}{c}
1,\left( \mathbf{x},t\right) \in G_{3h}, \\ 
0,\text{ }\left( \mathbf{x},t\right) \in \left( G_{h}\diagdown G_{2h}\right)
\cup \left( Q\left( R,T^{-},T\right) \diagdown G_{h}\right) , \\ 
\in \left[ 0,1\right] ,\text{ }\left( \mathbf{x},t\right) \in
G_{2h}\diagdown G_{3h}.%
\end{array}%
\right.%
\end{array}%
\right.  \label{5.18}
\end{equation}%
Such functions are constructed in the Analysis course. Multiply both sides
of equation (\ref{5.16}) by the function $\chi $ and denote 
\begin{equation}
v\left( \mathbf{x},t\right) =\chi \left( \mathbf{x},t\right) \widetilde{V}%
\left( \mathbf{x},t\right) .  \label{5.19}
\end{equation}%
We have 
\begin{equation}
\chi L_{0}\left( \widetilde{V}\right) =L_{0}\left( v\right) -2\chi _{t}%
\widetilde{V}_{t}-\chi _{tt}\widetilde{V}+2\nabla \chi \nabla \widetilde{V}%
+\Delta \chi \widetilde{V}.  \label{5.20}
\end{equation}%
By (\ref{5.16})-(\ref{5.20})%
\begin{equation}
\left. 
\begin{array}{c}
\left\vert L_{0}\left( v\right) \right\vert \leq C_{1}\left( \left\vert
\nabla v\right\vert +\left\vert v_{t}\right\vert +\left\vert v\right\vert
\right) +C_{1}\left( 1-\chi \right) \left( \left\vert \nabla \widetilde{V}%
\right\vert +\left\vert \widetilde{V}_{t}\right\vert +\left\vert \widetilde{V%
}\right\vert \right) \text{ in }G_{h}, \\ 
v\mid _{S\left( R,T^{-},T\right) }=\chi \widetilde{q}_{0},\text{ }\partial
_{\nu }v\mid _{S\left( R,T^{-},T\right) }=\chi \widetilde{q}_{1}+\partial
_{\nu }\chi \cdot \widetilde{q}_{0}.%
\end{array}%
\right.  \label{5.21}
\end{equation}

Square both sides of the inequality in the first line of (\ref{5.21}). Next,
recall (\ref{5.9}) and multiply both sides of the obtained inequality by the
function $\varphi ^{2}\left( \mathbf{x},t\right) .$ Next, using
Cauchy-Schwarz inequality and integrating over the domain $G_{h},$ we obtain%
\begin{equation}
\left. 
\begin{array}{c}
C_{1}e^{6\lambda h}\int\limits_{G_{h}\diagdown G_{3h}}\left( \left\vert
\nabla \widetilde{V}\right\vert ^{2}+\left\vert \widetilde{V}_{t}\right\vert
^{2}+\left\vert \widetilde{V}\right\vert ^{2}\right) d\mathbf{x}dt+ \\ 
+C_{1}\int\limits_{G_{h}}\left( \left\vert \nabla v\right\vert
^{2}+\left\vert v_{t}\right\vert ^{2}+\left\vert v\right\vert ^{2}\right)
\varphi ^{2}d\mathbf{x}dt\geq \\ 
\geq \int\limits_{G_{h}}\left\vert L_{0}\left( v\right) \right\vert
^{2}\varphi ^{2}d\mathbf{x}dt.%
\end{array}%
\right.  \label{5.22}
\end{equation}%
We now estimate from the below the term in the third line of (\ref{5.22}).
To do this, we integrate the Carleman estimate (\ref{5.10}) with condition (%
\ref{5.11}) over the domain $G_{h}.$ We apply Gauss formula to the terms $%
\text{div}U,W_{t}$ and take into account (\ref{5.4}), (\ref{5.81}), (\ref%
{5.12}), (\ref{5.18}), (\ref{5.19}) and (\ref{5.21}). In particular, $v\mid
_{\partial _{1}G_{h}}=\nabla v\mid _{\partial _{1}G_{h}}=v_{t}\mid
_{\partial _{1}G_{h}}=0.$ Hence, 
\begin{equation}
\left. 
\begin{array}{c}
\int\limits_{G_{h}}\left\vert L_{0}\left( v\right) \right\vert ^{2}\varphi
^{2}d\mathbf{x}dt\geq C\lambda \int\limits_{G_{h}}\left( \left\vert
v_{t}\right\vert ^{2}+\left\vert \nabla v\right\vert ^{2}+\lambda
^{2}\left\vert v\right\vert ^{2}\right) \varphi ^{2}d\mathbf{x}dt- \\ 
-C_{1}e^{2\lambda \left( 2R+\sigma \right) ^{2}}\left( \left\Vert \widetilde{%
q}_{0}\right\Vert _{H_{N}^{1}\left( S\left( R,T^{-},T\right) \right)
}^{2}+\left\Vert \widetilde{q}_{1}\right\Vert _{L_{2,N}\left( S\left(
R,T^{-},T\right) \right) }^{2}\right) ,\text{ }\forall \lambda \geq \lambda
_{0}.%
\end{array}%
\right.  \label{5.23}
\end{equation}%
To obtain the term $e^{2\lambda \left( 2R+\sigma \right) ^{2}}$ here, we
have used the fact that by (\ref{5.1}), (\ref{5.4}) and (\ref{5.9}) $\varphi
^{2}\left( \mathbf{x},t\right) \leq e^{2\lambda \left( 2R+\sigma \right)
^{2}}$ in $G_{h}.$ Combining (\ref{5.22}) with (\ref{5.23}), we obtain%
\begin{equation}
\left. 
\begin{array}{c}
C_{1}e^{6\lambda h}\int\limits_{G_{h}\diagdown G_{3h}}\left( \left\vert
\nabla \widetilde{V}\right\vert ^{2}+\left\vert \widetilde{V}_{t}\right\vert
^{2}+\left\vert \widetilde{V}\right\vert ^{2}\right) d\mathbf{x}dt+ \\ 
+C_{1}\int\limits_{G_{h}}\left( \left\vert \nabla v\right\vert
^{2}+\left\vert v_{t}\right\vert ^{2}+\left\vert v\right\vert ^{2}\right)
\varphi ^{2}d\mathbf{x}dt\geq \\ 
\geq C\lambda \int\limits_{G_{h}}\left( \left\vert v_{t}\right\vert
^{2}+\left\vert \nabla v\right\vert ^{2}+\lambda ^{2}\left\vert v\right\vert
^{2}\right) \varphi ^{2}d\mathbf{x}dt- \\ 
-C_{1}e^{2\lambda \left( 2R+\sigma \right) ^{2}}\left( \left\Vert \widetilde{%
q}_{0}\right\Vert _{H_{N}^{1}\left( S\left( R,T^{-},T\right) \right)
}^{2}+\left\Vert \widetilde{q}_{1}\right\Vert _{L_{2,N}\left( S\left(
R,T^{-},T\right) \right) }\right) ,\text{ }\forall \lambda \geq \lambda _{0}.%
\end{array}%
\right.  \label{5.24}
\end{equation}%
Choose a sufficiently large number $\lambda _{1}=\lambda _{1}\left( Q\left(
R,T^{\pm }\right) ,M_{1},M_{2},N,h\right) \geq \lambda _{0}\geq 1$ such that 
$C\lambda /2\geq C_{1}$ for all $\lambda \geq \lambda _{1}.$ Hence, (\ref%
{5.24}) implies%
\begin{equation}
\left. 
\begin{array}{c}
C_{1}e^{6\lambda h}\int\limits_{G_{h}\diagdown G_{3h}}\left( \left\vert
\nabla \widetilde{V}\right\vert ^{2}+\left\vert \widetilde{V}_{t}\right\vert
^{2}+\left\vert \widetilde{V}\right\vert ^{2}\right) d\mathbf{x}dt\geq \\ 
\geq \lambda \int\limits_{G_{h}}\left( \left\vert v_{t}\right\vert
^{2}+\left\vert \nabla v\right\vert ^{2}+\lambda ^{2}\left\vert v\right\vert
^{2}\right) \varphi ^{2}d\mathbf{x}dt- \\ 
-C_{1}e^{2\lambda \left( 2R+\sigma \right) ^{2}}\left( \left\Vert \widetilde{%
q}_{0}\right\Vert _{H_{N}^{1}\left( S\left( R,T^{-},T\right) \right)
}^{2}+\left\Vert \widetilde{q}_{1}\right\Vert _{L_{2,N}\left( S\left(
R,T^{-},T\right) \right) }^{2}\right) ,\text{ }\forall \lambda \geq \lambda
_{1}.%
\end{array}%
\right.  \label{5.25}
\end{equation}%
By (\ref{5.80}) $G_{4h}\subset G_{3h}\subset G_{h}$. Next, by (\ref{5.3}), (%
\ref{5.4}) and (\ref{5.9}) $\varphi \left( \mathbf{x},t\right) \geq
e^{8\lambda h}$ in $G_{4h}.$ Also, by (\ref{5.18}) and (\ref{5.19}) $v\left( 
\mathbf{x},t\right) =\widetilde{V}\left( \mathbf{x},t\right) $ in $G_{3h}.$
Hence, (\ref{5.25}) implies%
\begin{equation}
\left. 
\begin{array}{c}
C_{1}e^{2\lambda \left( 2R+\sigma \right) ^{2}}\left( \left\Vert \widetilde{q%
}_{0}\right\Vert _{H^{1}\left( S\left( R,T^{-},T\right) \right)
}^{2}+\left\Vert \widetilde{q}_{1}\right\Vert _{L_{2}\left( S\left(
R,T^{-},T\right) \right) }^{2}\right) + \\ 
+C_{1}e^{-2\lambda h}\left\Vert \widetilde{V}\right\Vert _{Q\left(
R,T^{-},T\right) }^{2}\geq \\ 
\geq \int\limits_{G_{4h}}\left( \left\vert \widetilde{V}_{t}\right\vert
^{2}+\left\vert \nabla \widetilde{V}\right\vert ^{2}+\left\vert \widetilde{V}%
\right\vert ^{2}\right) d\mathbf{x}dt,\text{ }\forall \lambda \geq \lambda
_{1}.%
\end{array}%
\right.  \label{5.26}
\end{equation}%
It follows from the mean value theorem, (\ref{5.80}) and (\ref{5.26}) that
there exists a number $t_{0}\in \left( T_{0}-A,T_{0}+A\right) $ such that%
\begin{equation}
\left. 
\begin{array}{c}
\left\Vert \widetilde{V}_{t}\left( x,t_{0}\right) \right\Vert
_{L_{2,N}\left( \Omega \right) }^{2}+\left\Vert \widetilde{V}\left(
x,t_{0}\right) \right\Vert _{H_{N}^{1}\left( \Omega \right) }^{2}\leq
C_{1}e^{-2\lambda h}\left\Vert \widetilde{V}\right\Vert _{H_{N}^{1}\left(
Q\left( R,T^{-},T\right) \right) }^{2}+ \\ 
+C_{1}e^{2\lambda \left( 2R+\sigma \right) ^{2}}\left( \left\Vert \widetilde{%
q}_{0}\right\Vert _{H_{N}^{1}\left( S\left( R,T^{-},T\right) \right)
}^{2}+\left\Vert \widetilde{q}_{1}\right\Vert _{L_{2,N}\left( S\left(
R,T^{-},T\right) \right) }^{2}\right) ,\text{ }\forall \lambda \geq \lambda
_{1}.%
\end{array}%
\right.  \label{5.27}
\end{equation}%
We recall now (\ref{1.4})-(\ref{1.50}) and apply the standard energy
estimate to the system of hyperbolic PDEs (\ref{5.16}) with the boundary
data (\ref{5.160}). Since the wave equation can be solved in both positive
and negative directions of time, then we apply the energy estimate twice: in
the domain $\Omega \times \left( t_{0},T\right) $ and in the domain $\Omega
\times \left( T^{-},t_{0}\right) .$ In both cases $\widetilde{V}\left(
x,t_{0}\right) $ and $\widetilde{V}_{t}\left( x,t_{0}\right) $ are used as
initial conditions. Hence, using (\ref{5.27}), we obtain%
\begin{equation}
\left. 
\begin{array}{c}
\left( 1-C_{1}e^{-\lambda h}\right) \left\Vert \widetilde{V}\right\Vert
_{H_{N}^{1}\left( Q\left( R,T^{-},T\right) \right) }\leq \\ 
\leq C_{1}e^{\lambda \left( 2R+\sigma \right) ^{2}}\left( \left\Vert 
\widetilde{q}_{0}\right\Vert _{H_{N}^{1}\left( S\left( R,T^{-},T\right)
\right) }+\left\Vert \widetilde{q}_{1}\right\Vert _{L_{2,N}\left( S\left(
R,T^{-},T\right) \right) }\right) ,\text{ }\forall \lambda \geq \lambda _{1}.%
\end{array}%
\right.  \label{5.28}
\end{equation}%
Fix a number $\lambda =\lambda _{2}\geq \lambda _{1}$ such that $%
1-C_{1}e^{-\lambda _{2}h}\geq 1/2.$ Then, for this value of $\lambda ,$ (\ref%
{5.28}) implies the target estimate (\ref{5.15}). $\square $

\section{Convexification}

\label{sec:5}

To solve BVP\ (\ref{4.5})-(\ref{4.7}) numerically, we use the
convexification method. We assume in this sections 5 that\emph{\ }conditions
(\ref{2.01})-(\ref{1.50}), (\ref{4.06}) and (\ref{4.9}) hold. Let $M_{1}$
and $M_{2}$ be two numbers defined in (\ref{4.8}) and (\ref{4.9}), and let $%
K\geq \max \left( M_{1},M_{2}\right) $ be an arbitrary number. Consider the
set $B\left( K\right) \subset H_{N}^{4}\left( Q\left( R,T^{-},T\right)
\right) $ of $N-$D vector functions,%
\begin{equation}
B\left( K\right) =\left\{ 
\begin{array}{c}
V\in H_{N}^{4}\left( Q\left( R,T^{-},T\right) \right) :\left\Vert
V\right\Vert _{H_{N}^{4}\left( Q\left( R,T^{-},T\right) \right) }<K, \\ 
V\mid _{S\left( R,T^{-},T\right) }=q_{0}\left( \mathbf{x},t\right) ,\text{ }%
\partial _{\nu }V\mid _{S\left( R,T^{-},T\right) }=q_{1}\left( \mathbf{x}%
,t\right)%
\end{array}%
\right\} .  \label{6.1}
\end{equation}%
By Sobolev embedding theorem%
\begin{equation}
B\left( K\right) \subset C_{N}^{1}\left( \overline{Q\left( R,T^{-},T\right) }%
\right) ,\left\Vert V\right\Vert _{C_{N}^{1}\left( \overline{Q\left(
R,T^{-},T\right) }\right) }\leq CK,\text{ }\forall V\in B\left( K\right) .
\label{6.2}
\end{equation}

Let $\varphi \left( \mathbf{x},t\right) $ and $\chi \left( \mathbf{x}%
,t\right) $ be the functions defined in (\ref{5.9}) and (\ref{5.18})
respectively.\ Denote 
\begin{equation}
F_{1}\left( \nabla V,V\right) =M_{N}^{-1}F\left( \nabla V,V\right) .
\label{6.3}
\end{equation}%
Consider the following functional%
\begin{equation}
\left. 
\begin{array}{c}
J_{\lambda ,\alpha }\left( V\right) =e^{-8\lambda h}\int\limits_{Q\left(
R,T^{-},T\right) }\left( \chi \left( V_{tt}-\Delta V\right) +\chi
F_{1}\left( \nabla V,V\right) \right) ^{2}\varphi ^{2}d\mathbf{x}dt+ \\ 
+\alpha \left\Vert V\right\Vert _{H_{N}^{4}\left( Q\left( R,T^{-},T\right)
\right) }^{2},%
\end{array}%
\right.  \label{6.4}
\end{equation}%
where $\alpha \in \left( 0,1\right) $ is the regularization parameter. The
multiplier $e^{-8\lambda h}$ balances first and second terms in the right
hand side of (\ref{6.4}), also, see (\ref{5.18}). We solve the following
Minimization Problem:

\textbf{Minimization Problem}. Minimize the functional $J_{\lambda ,\alpha
}\left( V\right) $ on the set $\overline{B\left( K\right) }.$

\subsection{The strong convexity of the functional $J_{\protect\lambda , 
\protect\alpha }\left( V\right) $}

\label{sec:5.1}

Denote $\left[ ,\right] $ the scalar product in the space $H_{N}^{4}\left(
Q\left( R,T^{-},T\right) \right) .$ Introduce the subspace

$H_{0,N}^{4}\left( Q\left( R,T^{-},T\right) \right) $ of the space $%
H_{N}^{4}\left( Q\left( R,T^{-},T\right) \right) $ as%
\begin{equation}
\left. H_{0,N}^{4}\left( Q\left( R,T^{-},T\right) \right) =\left\{ 
\begin{array}{c}
V\in H_{N}^{4}\left( Q\left( R,T^{-},T\right) \right) : \\ 
V\mid _{S\left( R,T^{-},T\right) }=\partial _{\nu }V\mid _{S\left(
R,T^{-},T\right) }=0%
\end{array}%
\right\} \right.  \label{6.5}
\end{equation}

\textbf{Theorem 5.1}. \emph{Assume that conditions (\ref{2.01})-(\ref{1.50}
), (\ref{4.06}), (\ref{4.9}) as well as conditions \ref{5.001}), (\ref{5.01}
), (\ref{5.2}), (\ref{5.3}), (\ref{5.4}), (\ref{5.6}) and (\ref{5.9}) hold.
Then:}

\emph{1. At each point }$V\in \overline{B\left( K\right) }$\emph{\ and for
all }$\lambda \geq 0,\alpha \in \left( 0,1\right) $\emph{\ there exists Fr 
\'{e}chet derivative }$J_{\lambda ,\alpha }^{\prime }\left( V\right) \in
H_{0,N}^{4}\left( Q\left( R,T^{-},T\right) \right) $\emph{\ of the
functional }$J_{\lambda ,\alpha }\left( V\right) .$\emph{\ Furthermore, this
derivative is Lipschitz continuous on }$\overline{B\left( K\right) }$\emph{.
In other words, there exists a number }$D>0$\emph{\ depending on the same
parameters from which }$J_{\lambda ,\alpha }\left( V\right) $\emph{\ depends
such that} 
\begin{equation}
\left. 
\begin{array}{c}
\left\Vert J_{\lambda ,\alpha }^{\prime }\left( V_{1}\right) -J_{\lambda
,\alpha }^{\prime }\left( V_{2}\right) \right\Vert _{H_{N}^{4}\left( Q\left(
R,T^{-},T\right) \right) }\leq \\ 
\leq D\left\Vert V_{1}-V_{2}\right\Vert _{H_{N}^{4}\left( Q\left(
R,T^{-},T\right) \right) },\text{ }\forall V_{1},V_{2}\in \overline{B\left(
K\right) }.%
\end{array}%
\right.  \label{6.6}
\end{equation}

\emph{2. Let }$\lambda _{0}\geq 1$\emph{\ be the large parameter of Theorem
4.1. There exists another sufficiently large number }$\lambda _{1}=\lambda
_{1}\left( Q\left( R,T^{-},T\right) ,M_{1},M_{2},N,K,h\right) \geq \lambda
_{0}$\emph{\ depending only on listed parameters such for all }$\lambda \geq
\lambda _{1}$\emph{\ and for all }$\alpha $ \emph{such that} 
\begin{equation}
\alpha \in \left[ 2e^{-\lambda h},1\right)  \label{6.7}
\end{equation}%
\emph{the functional }$J_{\lambda ,\alpha }\left( V\right) $\emph{\ is
strongly convex on the set }$\overline{B\left( K\right) }.$\emph{\ More
precisely, the following inequality holds:}%
\begin{equation}
\left. 
\begin{array}{c}
J_{\lambda ,\alpha }\left( V_{1}\right) -J_{\lambda ,\alpha }\left(
V_{2}\right) -J_{\lambda ,\alpha }^{\prime }\left( V_{2}\right) \left(
V_{1}-V_{2}\right) \geq \\ 
\geq C_{1}\left\Vert V_{1}-V_{2}\right\Vert _{H_{N}^{1}\left( G_{4h}\right)
}^{2}+\left( \alpha /2\right) \left\Vert V_{1}-V_{2}\right\Vert
_{H_{N}^{4}\left( Q\left( R,T^{-},T\right) \right) }^{2}, \\ 
\text{ }\forall V_{1},V_{2}\in \overline{B\left( K\right) },\text{ }\forall
\lambda \geq \lambda _{1}.%
\end{array}%
\right.  \label{6.8}
\end{equation}

\emph{3. For any value of }$\lambda \geq \lambda _{1}$\emph{\ there exists
unique minimizer }$V_{\min ,\lambda }\in $\emph{\ }$\overline{B\left(
K\right) }$\emph{\ of the functional }$J_{\lambda ,\alpha }\left( V\right) $%
\emph{\ on the set }$\overline{B\left( K\right) }.$\emph{\ Furthermore, the
following inequality holds:}%
\begin{equation}
\left[ J_{\lambda ,\alpha }^{\prime }\left( V_{\min ,\lambda }\right)
,V_{\min ,\lambda }-V\right] \leq 0,\text{ }\forall V\in \overline{B\left(
K\right) }.  \label{6.9}
\end{equation}

\textbf{Proof}. Let $V_{1},V_{2}\in \overline{B\left( K\right) }$ be two
arbitrary points. Denote $r=V_{1}-V_{2}.$ Then (\ref{6.1}), (\ref{6.5}) and
triangle inequality imply:%
\begin{equation}
\left. 
\begin{array}{c}
r\in \overline{B_{0}\left( 2K\right) }= \\ 
=\left\{ V:V\in H_{0,N}^{4}\left( Q\left( R,T^{\pm }\right) \right)
,\left\Vert V\right\Vert _{H_{0,N}^{4}\left( Q\left( R,T^{-},T\right)
\right) }\leq 2K\right\} .%
\end{array}%
\right.  \label{6.10}
\end{equation}
Since $V_{1}=V_{2}+r,$ then 
\begin{equation}
\left. 
\begin{array}{c}
\chi \left( \partial ^{2}V_{1}-\Delta V_{1}\right) +\chi F_{1}\left( \nabla
V_{1},V_{1}\right) = \\ 
=\left( \chi \left( \partial ^{2}V_{2}-\Delta V_{2}\right) +\chi \left(
r_{tt}-\Delta r\right) +\chi F_{1}\left( \nabla V_{2}+\nabla
r,V_{2}+r\right) \right) ^{2}.%
\end{array}%
\right.  \label{6.11}
\end{equation}%
Using the multidimensional analog of Taylor formula \cite{V}, we obtain
similarly with the right hand side of (\ref{5.16})%
\begin{equation*}
F_{1}\left( \nabla V_{2}+\nabla r,V_{2}+r\right) =F_{1}\left( \nabla
V_{2},V_{2}\right) +Y_{1}\left( \mathbf{x},t\right) \nabla r+Y_{2}\left( 
\mathbf{x},t\right) r,
\end{equation*}%
where $Y_{1}\left( \mathbf{x},t\right) $ is a $N\times 3N$ matrix, $%
Y_{2}\left( \mathbf{x},t\right) $ is a $N\times N$ matrix with their
elements belonging to the space $C\left( \overline{Q\left( R,T^{-},T\right) }%
\right) ,$ and 
\begin{equation}
\left\vert Y_{1}\left( \mathbf{x},t\right) \right\vert ,\left\vert
Y_{2}\left( \mathbf{x},t\right) \right\vert \leq C_{1},\left( \mathbf{x}%
,t\right) \in \overline{Q\left( R,T^{-},T\right) }.  \label{6.110}
\end{equation}%
Hence, (\ref{6.11}) can be rewritten as%
\begin{equation*}
\left. 
\begin{array}{c}
\chi \left( \partial ^{2}V_{1}-\Delta V_{1}\right) +\chi F_{1}\left( \nabla
V_{1},V_{1}\right) = \\ 
=\chi \left( \partial ^{2}V_{2}-\Delta V_{2}\right) +\chi F_{1}\left( \nabla
V_{2},V_{2}\right) + \\ 
+\chi \left( r_{tt}-\Delta r\right) +\chi Y_{1}\left( \mathbf{x},t\right)
\nabla r+\chi Y_{2}\left( \mathbf{x},t\right) r.%
\end{array}%
\right.
\end{equation*}%
Hence, 
\begin{equation*}
\left. 
\begin{array}{c}
\left( \chi \left( \partial ^{2}V_{1}-\Delta V_{1}\right) +\chi F_{1}\left(
\nabla V_{1},V_{1}\right) \right) ^{2}= \\ 
=\left( \chi \left( \partial ^{2}V_{2}-\Delta V_{2}\right) +\chi F_{1}\left(
\nabla V_{2},V_{2}\right) \right) ^{2}+ \\ 
+2\left[ \chi \left( \partial ^{2}V_{2}-\Delta V_{2}\right) +\chi
F_{1}\left( \nabla V_{2},V_{2}\right) \right] \times \\ 
\times \left[ \chi \left( r_{tt}-\Delta r\right) +\chi Y_{1}\left( \mathbf{x}%
,t\right) \nabla r+\chi Y_{2}\left( \mathbf{x},t\right) r\right] + \\ 
+\left( \chi \left( r_{tt}-\Delta r\right) +\chi Y_{1}\left( \mathbf{x}%
,t\right) \nabla r+\chi Y_{2}\left( \mathbf{x},t\right) r\right) ^{2}.%
\end{array}%
\right.
\end{equation*}%
It follows from this formula that 
\begin{equation}
J_{\lambda ,\alpha }\left( V_{2}+r\right) -J_{\lambda ,\alpha }\left(
V_{2}\right) =J_{\lambda ,\alpha ,\text{linear}}\left( r\right) +J_{\lambda
,\alpha ,\text{nonlinear}}\left( r\right) ,  \label{6.12}
\end{equation}%
where%
\begin{equation}
\left. 
\begin{array}{c}
J_{\lambda ,\alpha ,\text{linear}}\left( r\right) =e^{-8\lambda h}\times \\ 
\times \int\limits_{Q\left( R,T^{\pm }\right) }2\left[ \chi \left( \partial
^{2}V_{2}-\Delta V_{2}\right) +\chi F_{1}\left( \nabla V_{2},V_{2}\right) %
\right] \times \\ 
\times \left[ \chi \left( r_{tt}-\Delta r\right) +\chi Y_{1}\left( \mathbf{x}%
,t\right) \nabla r+\chi Y_{2}\left( \mathbf{x},t\right) r\right] \varphi
^{2}d\mathbf{x}dt+2\alpha \left[ V_{2},r\right] ,%
\end{array}%
\right.  \label{6.13}
\end{equation}%
\begin{equation}
\left. 
\begin{array}{c}
J_{\lambda ,\alpha ,\text{nonlinear}}\left( r\right) =e^{-8\lambda h}\times
\\ 
\times \int\limits_{Q\left( R,T^{-},T\right) }\left( \chi \left(
r_{tt}-\Delta r\right) +\chi Y_{1}\left( \mathbf{x},t\right) \nabla r+\chi
Y_{2}\left( \mathbf{x},t\right) r\right) ^{2}\varphi ^{2}d\mathbf{x}dt+ \\ 
+\alpha \left\Vert r\right\Vert _{H^{4}\left( Q\left( R,T^{-},T\right)
\right) }^{2}.%
\end{array}%
\right.  \label{6.14}
\end{equation}

By (\ref{6.5}), (\ref{6.10}) and (\ref{6.13}) $J_{\lambda ,\alpha ,\text{
linear}}\left( r\right) :H_{0}^{4}\left( Q\left( R,T^{-},T\right) \right)
\rightarrow \mathbb{R}$ is a bounded linear functional. Hence, by Riesz
theorem there exists an element $\widetilde{J}_{\lambda ,\alpha ,\text{linear%
}}\in H_{0}^{4}\left( Q\left( R,T^{-},T\right) \right) $ such that 
\begin{equation}
J_{\lambda ,\alpha ,\text{linear}}\left( r\right) =\left[ \widetilde{J}%
_{\lambda ,\alpha ,\text{linear}},r\right] ,\forall r\in H_{0,N}^{4}\left(
Q\left( R,T^{-},T\right) \right) .  \label{6.15}
\end{equation}%
Next, it follows from (\ref{6.12}) and (\ref{6.14}) that 
\begin{equation*}
\frac{J_{\lambda ,\alpha }\left( V_{2}+r\right) -J_{\lambda ,\alpha }\left(
V_{2}\right) -\left[ \widetilde{J}_{\lambda ,\alpha ,\text{linear}},r\right] 
}{\left\Vert r\right\Vert _{H_{N}^{4}\left( Q\left( R,T^{-},T\right) \right)
}}=O\left( \left\Vert r\right\Vert _{H_{N}^{4}\left( Q\left(
R,T^{-},T\right) \right) }\right) ,
\end{equation*}%
as $\left\Vert r\right\Vert _{H_{N}^{4}\left( Q\left( R,T^{-},T\right)
\right) }\rightarrow 0.$ Hence, $\widetilde{J}_{\lambda ,\alpha ,\text{linear%
}}\in H_{0,N}^{4}\left( Q\left( R,T^{-},T\right) \right) $ is the Fr\'{e}%
chet derivative of the functional $J_{\lambda ,\alpha }$ at the point $%
V_{2}, $ 
\begin{equation}
\widetilde{J}_{\lambda ,\alpha ,\text{linear}}=J_{\lambda ,\alpha }^{\prime
}\left( V_{2}\right) .  \label{6.16}
\end{equation}%
The proof of the Lipschitz continuity property (\ref{6.6}) is similar with
the proof of Theorem 3.1 of \cite{Bak} and is, therefore, omitted.

We now prove the strong convexity property (\ref{6.8}). It follows from (\ref%
{6.12}), (\ref{6.15}) and (\ref{6.16}) that 
\begin{equation}
J_{\lambda ,\alpha }\left( V_{2}+r\right) -J_{\lambda ,\alpha }\left(
V_{2}\right) -J_{\lambda ,\alpha }^{\prime }\left( V_{2}\right) \left(
r\right) =J_{\lambda ,\alpha ,\text{nonlinear}}\left( r\right) .
\label{6.17}
\end{equation}%
Hence, we need to estimate from the below $J_{\lambda ,\alpha ,\text{
nonlinear}}\left( r\right) .$ Applying Cauchy-

Schwarz inequality to (\ref{6.14}) and using (\ref{6.110}), we obtain 
\begin{equation}
\left. 
\begin{array}{c}
J_{\lambda ,\alpha ,\text{nonlinear}}\left( r\right) \geq e^{-8\lambda
h}\int\limits_{Q\left( R,T^{-},T\right) }\left( \chi \left( r_{tt}-\Delta
r\right) \right) ^{2}\varphi ^{2}d\mathbf{x}dt- \\ 
-C_{1}e^{-8\lambda h}\int\limits_{Q\left( R,T^{-},T\right) }\left( \left(
\chi \nabla r\right) ^{2}+\left( \chi r\right) ^{2}\right) \varphi d\mathbf{x%
}dt+\alpha \left\Vert r\right\Vert _{H^{3}\left( Q\left( R,T^{-},T\right)
\right) }^{2}.%
\end{array}%
\right.  \label{6.18}
\end{equation}%
Denote 
\begin{equation}
s\left( \mathbf{x},t\right) =\chi \left( \mathbf{x},t\right) r\left( \mathbf{%
\ x},t\right) .  \label{6.19}
\end{equation}%
Then, using (\ref{5.18}) and (\ref{6.19}) we obtain similarly with (\ref%
{5.20}) and (\ref{5.21})%
\begin{equation}
\left. 
\begin{array}{c}
J_{\lambda ,\alpha ,\text{nonlinear}}\left( r\right) \geq e^{-8\lambda
h}\int\limits_{G_{h}}\left( s_{tt}-\Delta s\right) ^{2}\varphi ^{2}d\mathbf{x%
}dt- \\ 
-C_{1}e^{-8\lambda h}\int\limits_{G_{h}}\left( \left( \nabla s\right)
^{2}+s^{2}\right) ^{2}\varphi ^{2}d\mathbf{x}dt- \\ 
-e^{-8\lambda h}\int\limits_{G_{h}\diagdown G_{3h}}\left( \left( \nabla
r\right) ^{2}+r^{2}\right) \varphi ^{2}d\mathbf{x}dt+\alpha \left\Vert
r\right\Vert _{H_{N}^{4}\left( Q\left( R,T^{-},T\right) \right) }^{2}\geq \\ 
\geq e^{-8\lambda h}\int\limits_{G_{h}}\left( s_{tt}-\Delta s\right)
^{2}\varphi ^{2}d\mathbf{x}dt- \\ 
-C_{1}e^{-8\lambda h}\int\limits_{G_{h}}\left( \left( \nabla s\right)
^{2}+s^{2}\right) \varphi ^{2}d\mathbf{x}dt- \\ 
-C_{1}e^{-2\lambda h}\left\Vert r\right\Vert _{H^{3}\left( Q\left(
R,T^{-},T\right) \right) }^{2}+\alpha \left\Vert r\right\Vert _{H^{3}\left(
Q\left( R,T^{-},T\right) \right) }^{2}.%
\end{array}%
\right.  \label{6.20}
\end{equation}%
Applying Carleman estimate of Theorem 4.1 to the first term in the third
line of (\ref{6.20}), we obtain%
\begin{equation*}
\left. 
\begin{array}{c}
J_{\lambda ,\alpha ,\text{nonlinear}}\left( r\right) \geq C\lambda
e^{-8\lambda h}\int\limits_{G_{h}}\left( \left( \nabla s\right)
^{2}+s_{t}^{2}+\lambda ^{2}s^{2}\right) \varphi ^{2}d\mathbf{x}dt- \\ 
-C_{1}e^{-8\lambda h}\int\limits_{G_{h}}\left( \left( \nabla s\right)
^{2}+s^{2}\right) \varphi ^{2}d\mathbf{x}dt- \\ 
-C_{1}e^{-2\lambda h}\left\Vert r\right\Vert _{H_{N}^{4}\left( Q\left(
R,T^{-},T\right) \right) }^{2}+\alpha \left\Vert r\right\Vert
_{H_{N}^{4}\left( Q\left( R,T^{-},T\right) \right) }^{2},\text{ }\forall
\lambda \geq \lambda _{0}\geq 1.%
\end{array}%
\right.
\end{equation*}%
Taking here $\lambda _{1}=\lambda _{1}\left( Q\left( R,T^{-},T\right)
,M_{1},M_{2},N,K,h\right) \geq \lambda _{0}$ so large that $C\lambda
_{1}>2C_{1},$ we obtain%
\begin{equation}
\left. 
\begin{array}{c}
J_{\lambda ,\alpha ,\text{nonlinear}}\left( r\right) \geq C_{1}\lambda
e^{-8\lambda h}\int\limits_{G_{h}}\left( \left( \nabla s\right)
^{2}+s_{t}^{2}+\lambda ^{2}s^{2}\right) \varphi ^{2}d\mathbf{x}dt- \\ 
-C_{1}e^{-2\lambda h}\left\Vert r\right\Vert _{H_{N}^{4}\left( Q\left(
R,T^{-},T\right) \right) }^{2}+\alpha \left\Vert r\right\Vert
_{H_{N}^{4}\left( Q\left( R,T^{-},T\right) \right) }^{2},\text{ }\forall
\lambda \geq \lambda _{1}.%
\end{array}%
\right.  \label{6.21}
\end{equation}%
Since by (\ref{5.80}) $Q\left( R,A\right) \subset G_{4h}\subset
G_{3h}\subset G_{h}$ and also since by (\ref{5.4}), (\ref{5.9}), (\ref{5.18}%
) and (\ref{6.19}) $\varphi ^{2}\left( \mathbf{x},t\right) \geq e^{8\lambda
h}$ in $G_{4h}$ as well as $s\left( \mathbf{x},t\right) =r\left( \mathbf{x}%
,t\right) $ in $G_{4h},$ then (\ref{6.21}) leads to%
\begin{equation}
\left. 
\begin{array}{c}
J_{\lambda ,\alpha ,\text{nonlinear}}\left( r\right) \geq
C_{1}\int\limits_{G_{4h}}\left( \left( \nabla r\right)
^{2}+r_{t}^{2}+r^{2}\right) d\mathbf{x}dt- \\ 
-C_{1}e^{-2\lambda h}\left\Vert r\right\Vert _{H_{N}^{4}\left( Q\left(
R,T^{-},T\right) \right) }^{2}+\alpha \left\Vert r\right\Vert
_{H_{N}^{4}\left( Q\left( R,T^{-},T\right) \right) }^{2},\text{ }\forall
\lambda \geq \lambda _{1}.%
\end{array}%
\right.  \label{6.22}
\end{equation}%
Using (\ref{6.7}), we obtain 
\begin{equation}
\left. 
\begin{array}{c}
-C_{1}e^{-2\lambda h}\left\Vert r\right\Vert _{H_{N}^{4}\left( Q\left(
R,T^{-},T\right) \right) }^{2}+\alpha \left\Vert r\right\Vert
_{H_{N}^{4}\left( Q\left( R,T^{-},T\right) \right) }^{2}= \\ 
=\left( \alpha /2\right) \left\Vert r\right\Vert _{H_{N}^{4}\left( Q\left(
R,T^{-},T\right) \right) }^{2}+\left( \alpha /2-C_{1}e^{-2\lambda h}\right)
\left\Vert r\right\Vert _{H_{N}^{4}\left( Q\left( R,T^{-},T\right) \right)
}^{2}\geq \\ 
\geq \left( \alpha /2\right) \left\Vert r\right\Vert _{H_{N}^{4}\left(
Q\left( R,T^{-},T\right) \right) }^{2}+\left( e^{-\lambda
h}-C_{1}e^{-2\lambda h}\right) \left\Vert r\right\Vert _{H_{N}^{4}\left(
Q\left( R,T^{-},T\right) \right) }^{2}\geq \\ 
\geq \left( \alpha /2\right) \left\Vert r\right\Vert _{H_{N}^{4}\left(
Q\left( R,T^{-},T\right) \right) }^{2},\text{ }\forall \lambda \geq \lambda
_{1}.%
\end{array}%
\right.  \label{6.23}
\end{equation}%
Combining (\ref{6.22}) with (\ref{6.23}) and using (\ref{6.17}), we obtain
the desired strong convexity estimate (\ref{6.8}).

As soon as (\ref{6.8}) is established, the existence and uniqueness of the
minimizer $V_{\min ,\lambda }$ as well as inequality (\ref{6.9}) follow
immediately from (\ref{6.6}) and a combination of Lemma 2.1 with Theorem 2.1
of \cite{Bak}. $\square $

\subsection{Estimating the distance between the minimizer $V_{\min ,\protect%
\lambda }$ and the exact solution}

\label{sec:5.2}

The input data for any inverse problem contain a noise. One of the main
concepts of the theory of Ill-Posed Problems is the assumption about the
existence of an exact solution of such a problem, i.e. a solution with
noiseless data \cite{T}. In terms of that theory, the above found minimizer $%
V_{\min ,\lambda }$ is called \textquotedblleft regularized solution". It is
of an obvious interest to estimate the distance between the regularized
solution and exact solution $V^{\ast }.$ To do this, one needs to introduce
a noise in the input data (\ref{4.7}), find the regularized solution for the
noisy data and then estimate the distance between this solution and $V^{\ast
}.$ Such estimates for noisy data were obtained for other versions of the
convexification method in a number of publications, see e.g. \cite%
{Bak,Klibgrad,KMFG,Kepid}. Furthermore, a version of such an estimate can be
obtained for problem (\ref{4.5})-(\ref{4.7}) as well. However, since we want
to simplify the presentation, we assume here that the data (\ref{4.7}) are
noiseless. The vector function $V_{\min ,\lambda }$ is still not the exact
solution then due to the presence of the regularization parameter $\alpha $
in the functional $J_{\lambda ,\alpha }.$ Thus, we estimate in this section
the distance between $V_{\min ,\lambda }$ and $V^{\ast }$ in the noiseless
case.

\textbf{Theorem 5.2.}\emph{\ Assume that there exists an exact solution }$%
V^{\ast }$\emph{\ of problem (\ref{4.5})-(\ref{4.7}), i.e. we assume that }$%
V^{\ast }$\emph{\ satisfies conditions (\ref{4.5})-(\ref{4.7}) and }$V^{\ast
}\in B\left( K\right) $\emph{. Assume that conditions (\ref{2.01})-(\ref%
{1.50}), (\ref{4.06}), (\ref{4.9}) as well as conditions (\ref{5.001})-(\ref%
{5.9}) hold. Let }$\lambda _{1}$\emph{\ be the value of }$\lambda $\emph{\
found in Theorem 5.1. For }$\lambda \geq \lambda _{1},$\emph{\ let }$V_{\min
,\lambda }\in \overline{B\left( K\right) }$\emph{\ be the minimizer of the
functional }$J_{\lambda ,\alpha }$\emph{\ on the set }$\overline{B\left(
K\right) },$\emph{\ which was found in Theorem 5.1. By (\ref{6.7}) let }$%
\alpha =2e^{-\lambda h}$\emph{. Then the following accuracy estimate holds:}%
\begin{equation}
\left\Vert V^{\ast }-V_{\min ,\lambda }\right\Vert _{H_{N}^{1}\left(
G_{4h}\right) }\leq C_{1}e^{-\lambda h}.  \label{6.230}
\end{equation}

\textbf{Proof.} Since both vector functions $V_{\min ,\lambda },V^{\ast }\in 
\overline{B\left( K\right) },$ then, using (\ref{6.8}), we obtain%
\begin{equation}
\left. 
\begin{array}{c}
J_{\lambda ,\alpha }\left( V^{\ast }\right) -J_{\lambda ,\alpha }\left(
V_{\min ,\lambda }\right) -J_{\lambda ,\alpha }^{\prime }\left( V_{\min
,\lambda }\right) \left( V^{\ast }-V_{\min ,\lambda }\right) \geq \\ 
\geq C_{1}\left\Vert V^{\ast }-V_{\min ,\lambda }\right\Vert
_{H_{N}^{1}\left( G_{4h}\right) }^{2}.%
\end{array}%
\right.  \label{6.24}
\end{equation}%
By (\ref{6.4}) and (\ref{6.9})%
\begin{equation*}
-J_{\lambda ,\alpha }\left( V_{\min ,\lambda }\right) -J_{\lambda ,\alpha
}^{\prime }\left( V_{\min ,\lambda }\right) \left( V^{\ast }-V_{\min
,\lambda }\right) \leq 0.
\end{equation*}%
Hence, (\ref{6.24}) implies%
\begin{equation}
\left\Vert V^{\ast }-V_{\min ,\lambda }\right\Vert _{H_{N}^{1}\left(
G_{4h}\right) }^{2}\leq J_{\lambda ,\alpha }\left( V^{\ast }\right) .
\label{6.25}
\end{equation}%
Since $\chi \left( V_{tt}^{\ast }-\Delta V^{\ast }\right) +\chi F_{1}\left(
\nabla V^{\ast },V^{\ast }\right) =0,$ then (\ref{6.4}) implies that

$J_{\lambda ,\alpha }\left( V^{\ast }\right) =\alpha \left\Vert V^{\ast
}\right\Vert _{H^{3}\left( Q\left( R,T^{\pm }\right) \right) }^{2}\leq
C_{1}\alpha .$ Next, since $\alpha =2e^{-\lambda h},$ then (\ref{6.25})
implies the target estimate (\ref{6.230}). $\ \square $

\subsection{Global convergence of the gradient descent method}

\label{sec:5.3}

Let conditions of Theorem 5.2 hold. For $\lambda \geq \lambda _{1},$ assume
that 
\begin{equation}
V^{\ast }\in B\left( \frac{K}{3}-C_{1}e^{-\lambda h}\right) .  \label{6.26}
\end{equation}%
Then (\ref{6.230}) implies that it is reasonable to assume that 
\begin{equation}
V_{\min ,\lambda }\in B\left( \frac{K}{3}\right) .  \label{6.27}
\end{equation}%
Let 
\begin{equation}
V_{0}\in B\left( \frac{K}{3}\right)  \label{6.28}
\end{equation}%
For a number $\gamma >0$ consider the gradient descent method of the
minimization of the functional $J_{\lambda ,\alpha },$%
\begin{equation}
V_{n}=V_{n-1}-\gamma J_{\lambda ,\alpha }^{\prime }\left( V_{n-1}\right) ,%
\text{ }n=1,2,...  \label{6.29}
\end{equation}%
Since by Theorem 5.1 $J_{\lambda ,\alpha }^{\prime }\left( V_{n-1}\right)
\in H_{0,N}^{4}\left( Q\left( R,T^{-},T\right) \right) ,$ then boundary
conditions (\ref{4.7}) are the same for all vector functions $V_{n}.$

\textbf{Theorem 5.3.} \emph{Assume that conditions of Theorem 5.2 as well as
conditions (\ref{6.26})-(\ref{6.29}) hold. Then there exists a sufficiently
small number }$\gamma _{0}\in \left( 0,1\right) $\emph{\ such that for any }$%
\gamma \in \left( 0,\gamma _{0}\right) $\emph{\ all vector functions }%
\begin{equation}
V_{n}\in B\left( K\right) .  \label{6.30}
\end{equation}%
\emph{\ Furthermore, the sequence }$\left\{ V_{n}\right\} _{n=1}^{\infty }$%
\emph{\ converges to the minimizer }$V_{\min ,\lambda }$\emph{\ and there
exists a number }$\theta \in \left( 0,1\right) $\emph{\ such that the
following convergence estimates are valid: }%
\begin{equation}
\left\Vert V_{n}-V_{\min ,\lambda }\right\Vert _{H_{N}^{4}\left( Q\left(
R,T^{-},T\right) \right) }\leq \theta ^{n}\left\Vert V_{0}-V_{\min ,\lambda
}\right\Vert _{H_{N}^{4}\left( Q\left( R,T^{-},T\right) \right) },
\label{6.31}
\end{equation}%
\begin{equation}
\left\Vert V^{\ast }-V_{n}\right\Vert _{H_{N}^{1}\left( G_{4h}\right) }\leq
C_{1}e^{-\lambda h}+\theta ^{n}\left\Vert V_{0}-V_{\min ,\lambda
}\right\Vert _{H_{N}^{4}\left( Q\left( R,T^{-},T\right) \right) }.
\label{6.32}
\end{equation}

\textbf{Proof. }Formulas (\ref{6.30}) and (\ref{6.31}) follow immediately
from Theorem \ 6 of\emph{\ }\cite{Klibgrad}. Using (\ref{6.230}), (\ref{6.31}%
) and triangle inequality, we now prove (\ref{6.32}),%
\begin{equation*}
\left. 
\begin{array}{c}
\left\Vert V^{\ast }-V_{n}\right\Vert _{H_{N}^{1}\left( G_{4h}\right) }\leq
\left\Vert V^{\ast }-V_{\min ,\lambda }\right\Vert _{H_{N}^{1}\left(
G_{4h}\right) }+\left\Vert V_{n}-V_{\min ,\lambda }\right\Vert
_{H_{N}^{1}\left( G_{4h}\right) }\leq \\ 
\leq C_{1}e^{-\lambda h}+\left\Vert V_{n}-V_{\min ,\lambda }\right\Vert
_{H_{N}^{4}\left( Q\left( R,T^{-},T\right) \right) }\leq \\ 
\leq C_{1}e^{-\lambda h}+\theta ^{n}\left\Vert V_{0}-V_{\min ,\lambda
}\right\Vert _{H^{3}\left( Q\left( R,T^{-},T\right) \right) }.%
\end{array}%
\right.
\end{equation*}%
$\square $

\textbf{Remarks 5.1:}

\begin{enumerate}
\item Since the smallness of the number $K$ is not assumed in Theorem 5.3,
then this theorem claims the global convergence of the gradient descent
method (\ref{6.29}) of the minimization of the functional $J_{\lambda
,\alpha },$ see Definition in section 1.

\item Even though our convergence theory is constructed only for
sufficiently large values of the parameter $\lambda ,$ our computational
experience of this and past publications about the convexification method
demonstrates that actually reasonable values of the parameter $\lambda \in %
\left[ 1,5\right] $ work well, see the next section and \cite%
{Bak,KL,Klibgrad,KMFG,Kepid}. This is actually similar with many asymptotic
theories. Indeed, any such theory claims that if a certain parameter $X$ is
sufficiently large, then a certain formula $Y$ is accurate. However, in
practical computations only numerical experiments can define which exactly
values of $X$ provide a good accuracy of $Y$.
\end{enumerate}

\section{Numerical Studies}

\label{sec:6}

For the ball $\Omega $ in (\ref{2.01}), we choose $R=1/2$. We choose in (\ref%
{2.04}) and (\ref{1.40}) 
\begin{equation*}
T=12>T^{-}=4>R(\sqrt{5}+1)=\frac{\sqrt{5}+1}{2},
\end{equation*}%
\begin{equation*}
A=\frac{4}{3}.
\end{equation*}%
We set in (\ref{5.6}) 
\begin{equation*}
\sigma =2.5>\frac{2R}{\sqrt{2}-1}=\frac{1}{\sqrt{2}-1}.
\end{equation*}
In (\ref{5.01}) we choose 
\begin{equation*}
h=0.1\in \left( 0,\frac{\sigma ^{2}}{5}\right) =\left( 0,1.25\right) .
\end{equation*}
Then in (\ref{5.2})%
\begin{equation*}
\eta =\frac{\sigma ^{2}-h}{\left( 2A\right) ^{2}}=\frac{1107}{1280}\in
\left( 0,1\right) .
\end{equation*}

For function $a\left( \mathbf{x},t\right) $ that characterizes the moving
target in (\ref{1.9}), we set 
\begin{equation}
a\left( \mathbf{x},t\right)=\left\{ 
\begin{array}{cc}
a_{0}=const.>1, & \mbox{inside the moving target,} \\ 
1, & \mbox{outside the moving target.}%
\end{array}
\right.  \label{6.01}
\end{equation}
And we define 
\begin{align}
\mbox{correct inclusion/background contrast} &=\frac{a_{0}}{1},  \label{6.02}
\\
\mbox{computed inclusion/background contrast} &=\frac{\max_{ \mbox{inclusion}%
}\left( a _{\mbox{comp}}(\mathbf{x}, t)\right) }{1},  \label{6.03}
\end{align}%
where $a _{\mbox{comp}}(\mathbf{x}, t)$ is the reconstructed coefficient $a(%
\mathbf{x}, t).$ In the numerical tests below, we take $a_{0}=2$, which
corresponds to 2:1 inclusion/background contrasts. We consider three cases
of the moving target.

We solve the forward problem (\ref{1.12})-(\ref{1.13}) numerically by the
finite difference method with spatial mesh sizes $h_{x}=1/160\times 1/160$
and temporal mesh step size $h_{t}=1/640$. To generate the discrete point
source running along the interval $L_{0}$ in (\ref{2.03}), we choose $%
s=-R+ih_{s},i=1,2,3,\cdots ,100$ with $h_{s}=2R/101=1/101$.

To minimize the functional in (\ref{6.4}), we choose the spatial mesh sizes $%
\widetilde{h}_{x}=1/20\times 1/20$ and temporal mesh step size $\widetilde{h}%
_{t}=1/10$. We use the finite differences and numerical integration to
approximate the differential operators and integration in . Then we minimize
the functional $J_{\lambda ,\alpha }\left( V\right) $ with respect to the
values of the vector function $V\left( \mathbf{x},t\right) $ at the grid
points. We choose 
\begin{equation*}
N=5,\text{ }\lambda =3\text{ and }\alpha =0.01.
\end{equation*}%
To simplify our computational goal, we set $\chi \left( \mathbf{x},t\right)
\equiv 1$ in (\ref{6.4}) in our computations. In addition, we replace in (%
\ref{6.4}) $H_{N}^{4}\left( Q\left( R,T^{-},T\right) \right) $ with $%
H_{N}^{2}\left( Q\left( R,T^{-},T\right) \right) $ in our computations. The
fact that, regardless on these replacements, we still obtain accurate
numerical results, manifests a good degree of robustness of the numerical
method of this paper.

Tests 1 and 2 below are computed for the noiseless input data. However, in
the mos challenging Test 3 we introduce the random noise in the observation
data in \eqref{1.190} as follows: 
\begin{equation}
\left. 
\begin{array}{c}
g_{0}^{\xi _{0}}\left( \mathbf{x},\mathbf{x}_{0},t\right) =g_{0}\left( 
\mathbf{x},\mathbf{x}_{0},t\right) \left( 1+\delta \xi _{0}\left( \mathbf{x}%
,t\right) \right) ,\quad \left( \mathbf{x},t\right) \in S\left(
R,T^{-},T\right) ,\quad \mathbf{x}_{0}\in L_{0}, \\ 
g_{1}^{\xi _{1}}\left( \mathbf{x},\mathbf{x}_{0},t\right) =g_{1}\left( 
\mathbf{x},\mathbf{x}_{0},t\right) \left( 1+\delta \xi _{1}\left( \mathbf{x}%
,t\right) \right) ,\quad \left( \mathbf{x},t\right) \in S\left(
R,T^{-},T\right) ,\quad \mathbf{x}_{0}\in L_{0},%
\end{array}%
\right.  \label{6.1001}
\end{equation}%
where $\xi _{0}$ and $\xi _{1}$ are the uniformly distributed random
variables in the interval $[-1,1]$ depending on the value $\left( \mathbf{x}%
,t\right) \in S\left( R,T^{-},T\right) $. We choose $\delta =0.03,$ which
corresponds to the $3\%$ noise level. To calculate the first $s-$derivatives
of the noisy functions $p_{0}\left( \mathbf{x},\mathbf{x}_{0},t\right) $ and 
$p_{1}\left( \mathbf{x},\mathbf{x}_{0},t\right) $, we use the natural cubic
splines in $(-R,R)$ with the mesh grid size $h_{s}=R/100$ to approximate the
noisy input data (\ref{6.1001}), and then use the derivatives of those
splines to approximate the derivatives of corresponding noisy observation
data.

To guarantee that the solution of the problem of the minimization of
functional (\ref{6.4}) satisfies the boundary conditions (\ref{4.7}), we
adopt the Matlab's built-in optimization toolbox \textbf{fmincon} to
minimize the discretized form of these corresponding functions. The
iterations of \textbf{fmincon} stop when the followin condition is met: 
\begin{equation*}
|\nabla J_{\lambda ,\alpha }\left( V\right) |<10^{-2}.
\end{equation*}

\begin{figure}[tbph]
\centering
\includegraphics[width = 4.5in]{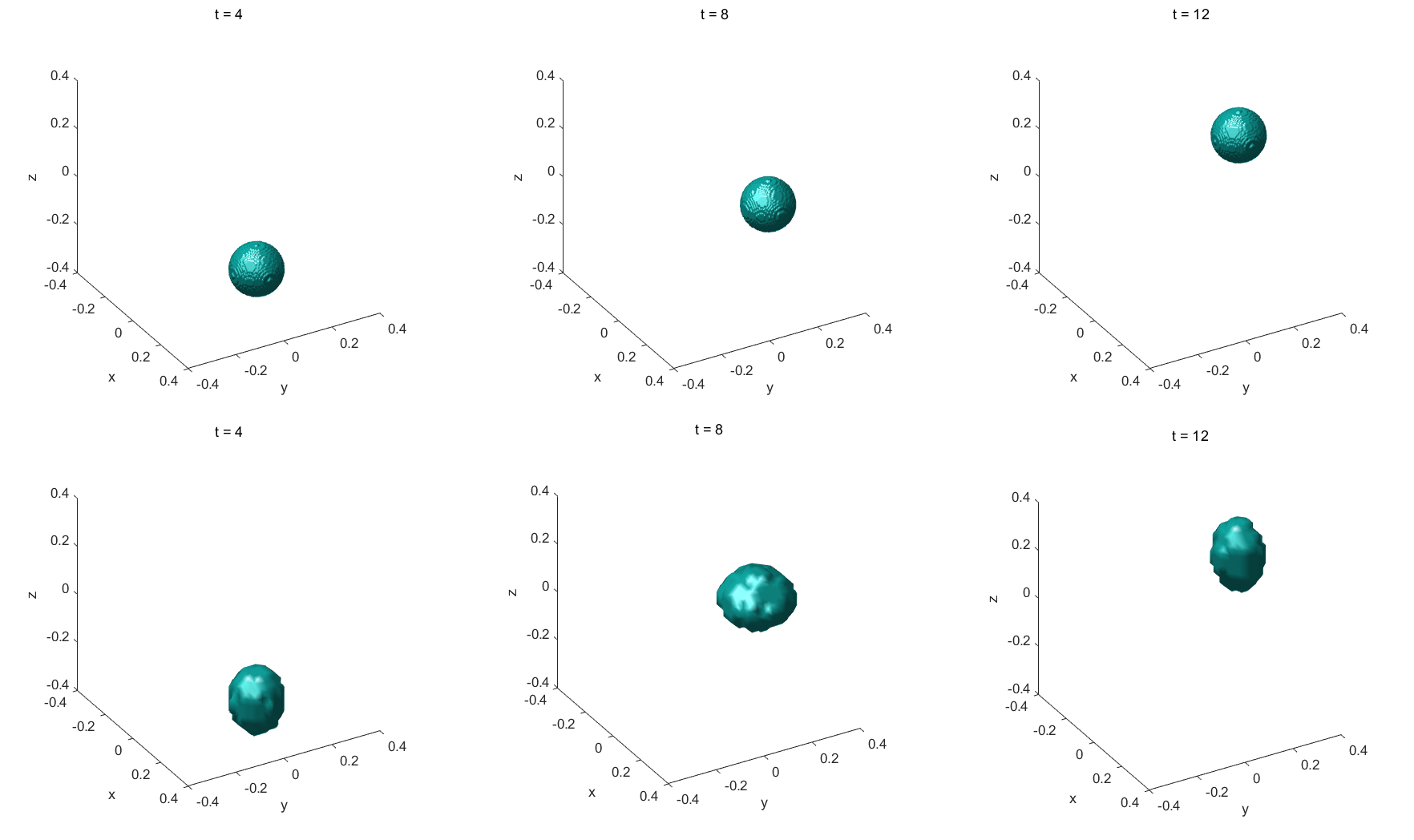}
\caption{Test 1. Exact (top) and reconstructed (bottom) moving ball, whose
center is given in (\protect\ref{6.04}).}
\label{ball_re}
\end{figure}

\begin{figure}[tbph]
\centering
\includegraphics[width = 4.5in]{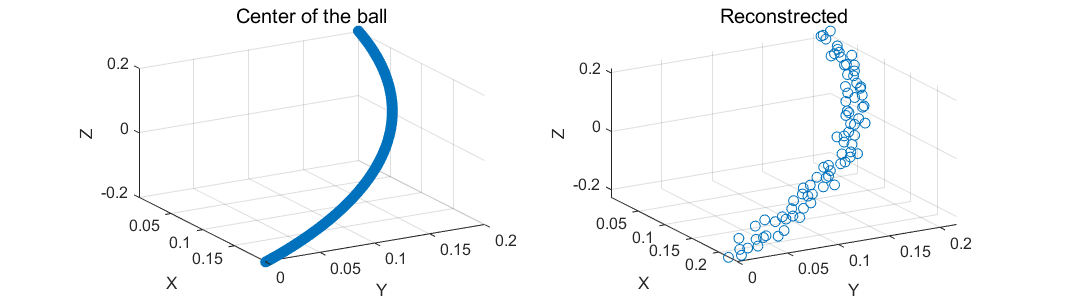}
\caption{Test 1. Exact (left) and reconstructed (right) center of the moving
ball, whose center is given in (\protect\ref{6.04}).}
\label{ball_center_re}
\end{figure}

\textbf{Test 1.} For the first moving target, we consider a moving ball,
whose radius is 0.1 and the center has the following trajectory:%
\begin{equation}
\left. 
\begin{array}{c}
x(t)=0.2\cos \left( (t-4)\pi /16\right) ,\quad y(t)=0.2\sin ((t-4)\pi /16),
\\ 
z(t)=0.05t-0.4,\quad t\in \left( 4,12\right) .%
\end{array}%
\right.  \label{6.04}
\end{equation}%
Both exact and reconstructed images are shown on Figure \ref{ball_re}. Exact
and reconstructed centers are shown on Figure \ref{ball_center_re}.

\textbf{Test 2.} For the second moving target, we consider a moving circular
cylinder, whose radius is 0.1 and the height is 0.2. The center of this
cylinder has the following trajectory: 
\begin{equation}
\left. 
\begin{array}{c}
x(t)=0.05t-0.4, \quad y(t)=0.05t-0.4, \quad z(t)=0.05t-0.4, \quad t\in
\left( 4, 12\right) .%
\end{array}%
\right.  \label{6.05}
\end{equation}

\begin{figure}[tbph]
\centering
\includegraphics[width = 4.5in]{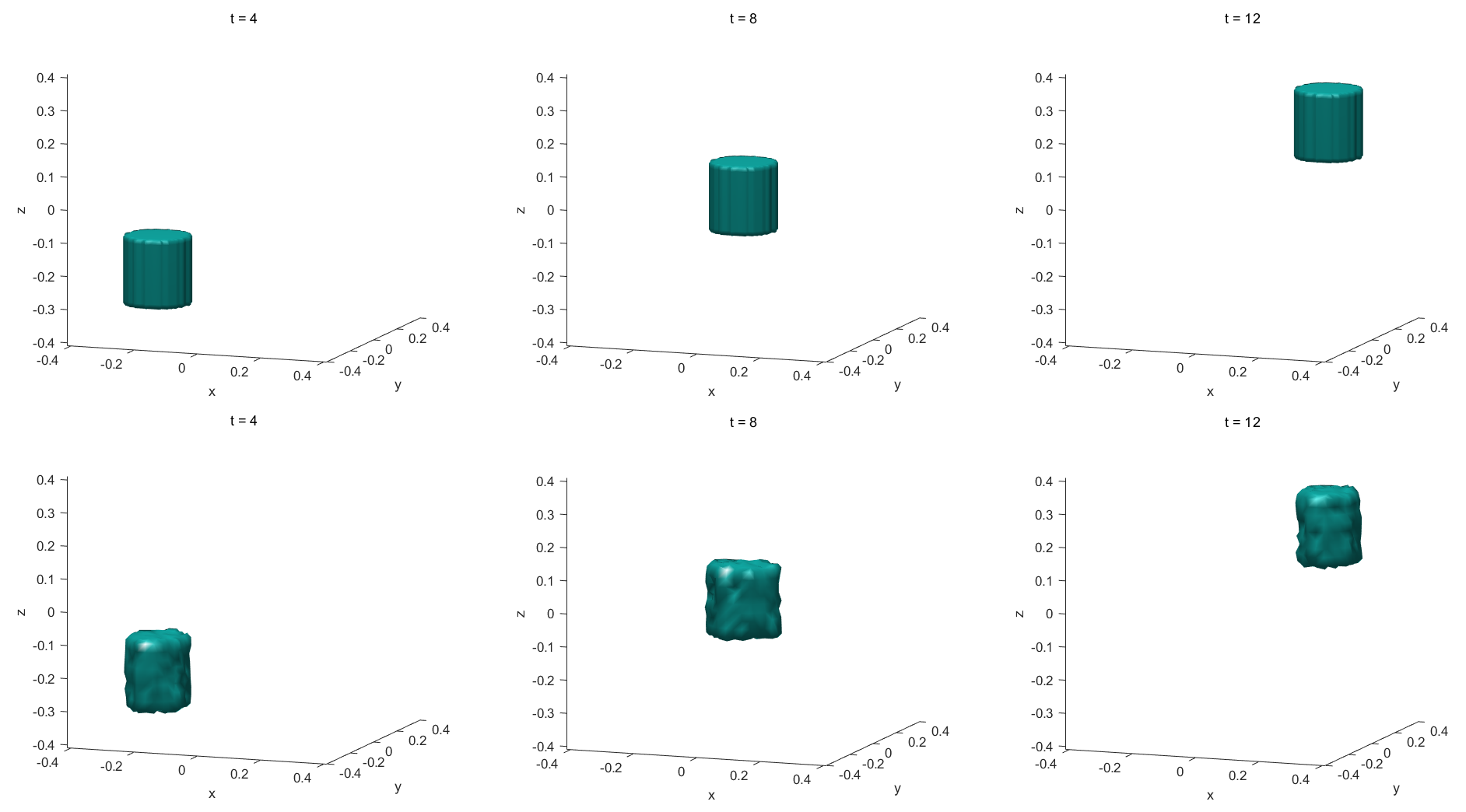}
\caption{Test 2. Exact (top) and reconstructed (bottom) moving cylinder,
whose center is given in (\protect\ref{6.05}).}
\label{cylinder_re}
\end{figure}

\begin{figure}[tbph]
\centering
\includegraphics[width = 4.5in]{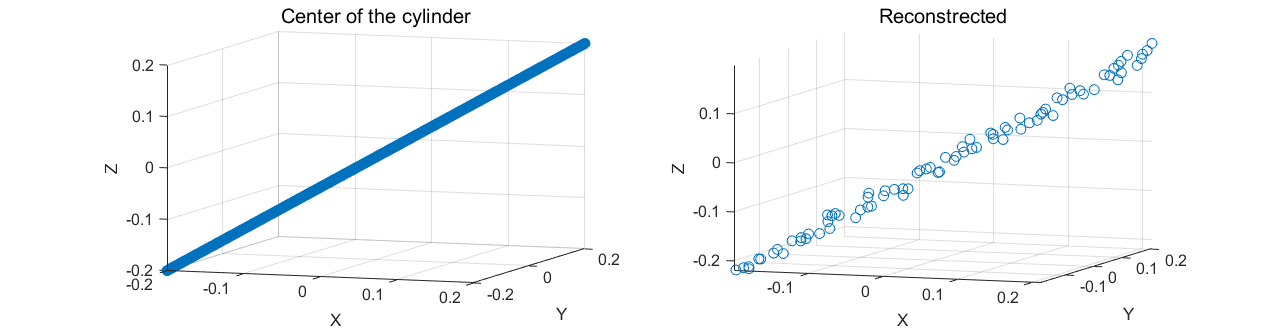}
\caption{Test 2. Exact (left) and reconstructed (right) center of the moving
cylinder, whose center is given in (\protect\ref{6.05}).}
\label{cylinder_center_re}
\end{figure}
Both exact and reconstructed images are shown on Figure \ref{cylinder_re}.
Exact and reconstructed centers are shown on Figure \ref{cylinder_center_re}.

\textbf{Test 3.} For the third moving target, we consider a moving circular
cylinder rotated along the $y-$axis with thr angle $\phi _{c}$. The radius
of this cylinder is 0.1 and the height is also 0.1. The center of this
cylinder has the following trajectory: 
\begin{equation}
\left. 
\begin{array}{c}
x(t)=0.4\sin ((t-4)\pi /16)-0.2,\quad y(t)=0.05t-0.4, \\ 
z(t)=-0.4\cos ((t-4)\pi /16)+0.2,\quad t\in \left( 4,12\right) .%
\end{array}%
\right.  \label{6.06}
\end{equation}%
The rotation angle $\phi _{c}$ is given as 
\begin{equation}
\phi _{c}(t)=t\pi /48+\pi /12,\quad t\in \left( 4,12\right) .  \label{6.07}
\end{equation}

In this test, we consider $\delta =0.03$ in the noisy observation data in %
\eqref{6.1001}.

\begin{figure}[tbph]
\centering
\includegraphics[width = 4.5in]{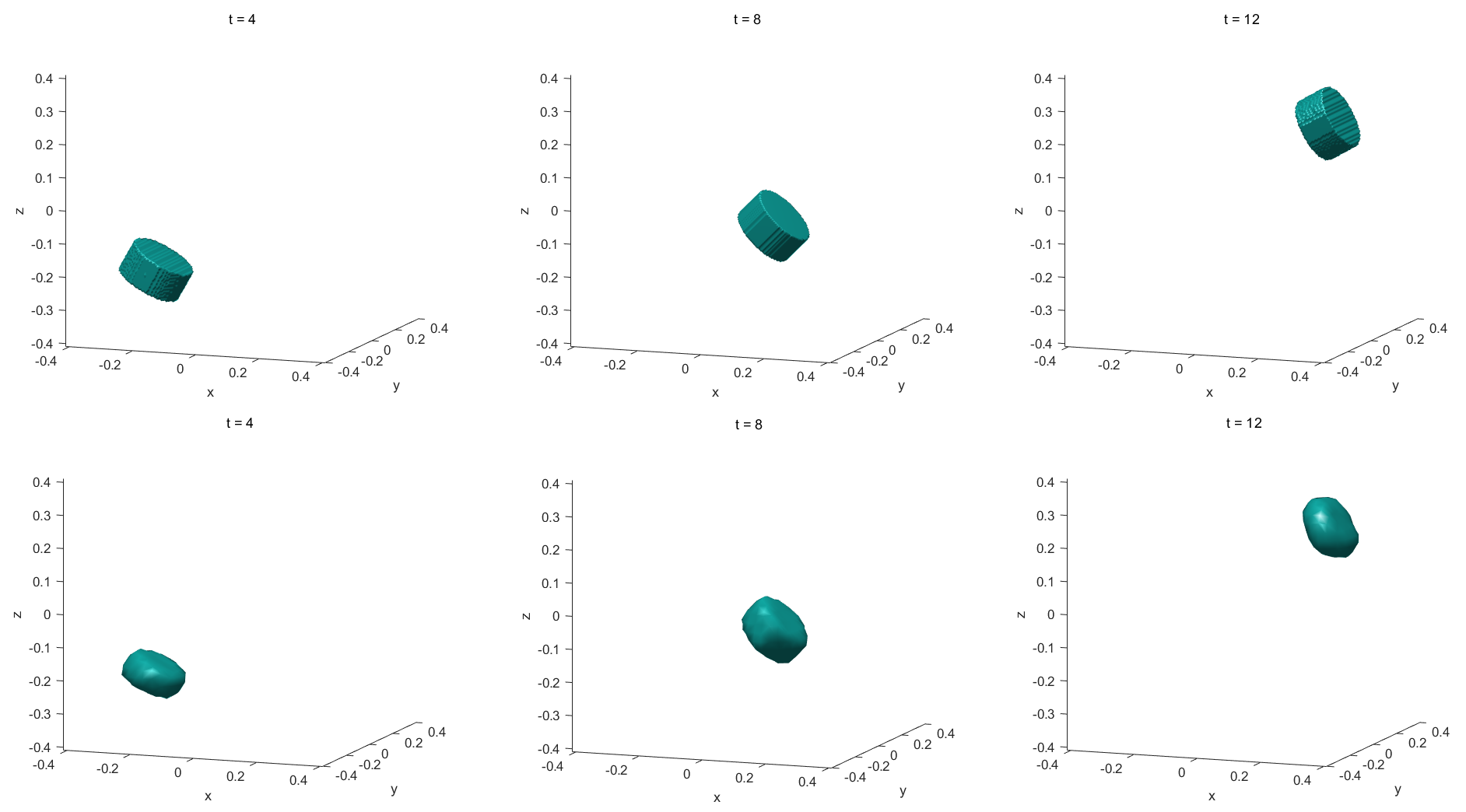}
\caption{Test 3. Exact (top) and reconstructed (bottom) moving rotated
cylinder, whose center is given in (\protect\ref{6.06}) and the rotation angle along $y-$axis is given in (\ref{6.07}).}
\label{rotated_cylinder_re}
\end{figure}

\begin{figure}[tbph]
\centering
\includegraphics[width = 4.5in]{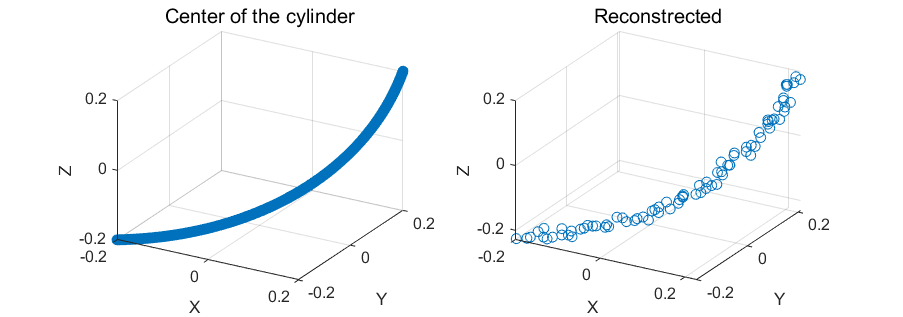}
\caption{Test 2. Exact (left) and reconstructed (right) center of the moving
cylinder, whose center is given in (\protect\ref{6.06}).}
\label{rotated_cylinder_center_re}
\end{figure}

Both exact and reconstructed images are shown on Figure \ref%
{rotated_cylinder_re}. Exact and reconstructed centers are shown on Figure %
\ref{rotated_cylinder_center_re}.

\textbf{Acknowledgments}. 
The work of Klibanov was partially supported by the US National Science Foundation grant DMS 2436227. 
The work of Li was partially supported by the Shenzhen Sci-Tech Fund No. RCJC20200714114556020, Guangdong Basic and Applied Research Fund No. 2023B1515250005. 
The work of Romanov was supported by the Mathematical Center in Akademgorodok under Agreement 075-15-2025-348 with the Ministry of Science and Higher Education of the Russian Federation. 
The work of Yang was partially supported by NSFC grant 12401558, Tianyuan Fund for Mathematics of the National Natural Science Foundation of China 12426105, Open Research Fund of Key Laboratory of Nonlinear Analysis \& Applications (Central China Normal University), Ministry of Education, P. R. China, and Supercomputing Center of Lanzhou University.


\section{Appendix: Proof of Theorem 2.1}

\label{sec:7}

By (\ref{1.16}) the problem (\ref{1.12}), (\ref{1.13}) in the domain $K_{,%
\mathbf{x}_{0},T}^{+}$ is equivalent to the following integral equation%
\begin{equation*}
u(\mathbf{x},\mathbf{x}_{0},t)=\frac{1}{4\pi -|\mathbf{x}-\mathbf{x}_{0}|}+
\end{equation*}%
\begin{equation}
+\frac{1}{4\pi }\int\limits_{D(\mathbf{x},\mathbf{x}_{0},t)}\frac{a(\xi
,t-|\xi -\mathbf{x}|)u(\xi ,\mathbf{x}_{0},t-|\xi -\mathbf{x}|)}{|\xi -%
\mathbf{x}|}\,d\xi ,~(\mathbf{x},t)\in K_{\mathbf{x}_{0},T}^{+},  \label{104}
\end{equation}%
where $D(\mathbf{x},\mathbf{x}_{0},t)$ is the inner part of the ellipsoid 
\begin{equation*}
E(\mathbf{x},\mathbf{x}_{0},t)=\{\xi \in \mathbb{R}^{3}|\,|\xi -\mathbf{x}%
|=t-|\xi -\mathbf{x}_{0}\}.
\end{equation*}%
Introduce curvilinear coordinates. To do this, we follow the work \cite{Rom}%
. Consider the family of ellipsoids 
\begin{equation*}
E(\mathbf{x},\mathbf{x}_{0},\tau )=\{\xi \in \mathbb{R}^{3}|\,|\xi -\mathbf{x%
}|+|\xi -\mathbf{x}_{0}|=\tau \},~\tau \in \big[|\mathbf{x}-\mathbf{x}_{0}|,t%
\big],
\end{equation*}%
and take the parameter $\tau $ as one of the curvilinear coordinates. Note
that if $\tau \rightarrow |\mathbf{x}-\mathbf{x}_{0}|+0,$ then $E(\mathbf{x},%
\mathbf{x}_{0},\tau )$ becomes the segment $L(\mathbf{x},\mathbf{x}_{0})$
joining points $\mathbf{x}$ and $\mathbf{x}_{0}$ . Let 
\begin{equation*}
\left( r,\varphi ,\theta \right) ,r>0,\varphi \in \lbrack 0,2\pi ),\theta
\in \lbrack 0,\pi ]
\end{equation*}%
be the spherical coordinate system with the center at the point $\mathbf{x}%
_{0}.$ Introduce the new coordinate system $\xi ^{\prime }=(\xi _{1}^{\prime
},\xi _{2}^{\prime },\xi _{3}^{\prime })$ with its center at $\mathbf{x}%
_{0}. $ Set%
\begin{equation}
\left. 
\begin{array}{c}
\xi =\mathbf{x}_{0}+\sum_{i=1}^{3}\xi _{i}^{\prime }\mathbf{e}_{i}, \\ 
\xi _{1}^{\prime }=r\sin \theta \cos \varphi ,~\xi _{2}^{\prime }=r\sin
\theta \sin \varphi ,~\xi _{3}^{\prime }=r\cos \theta ,%
\end{array}%
\right.  \label{105}
\end{equation}
where unit vectors $\mathbf{e}_{i},i=1,2,3,$ are chosen such that $\mathbf{e}%
_{3}=(\mathbf{x}-\mathbf{x}_{0})/|\mathbf{x}-\mathbf{x}_{0}|$ and all three
vectors $\mathbf{e}_{1},\mathbf{e}_{2},\mathbf{e}_{3}$ are orthogonal to
each other.

It follows from (\ref{105}) that 
\begin{equation*}
|\xi -\mathbf{x}_{0}|=r,~|\xi -\mathbf{x}|=(|\mathbf{x}-\mathbf{x}
_{0}|^{2}-2r|\mathbf{x}-\mathbf{x}_{0}|\cos \theta +r^{2})^{1/2}.
\end{equation*}%
Then the equation for the ellipsoid $E(\mathbf{x},\mathbf{x}_{0},\tau )$ is: 
\begin{equation}
r=\frac{\tau ^{2}-|\mathbf{x}-\mathbf{x}_{0}|^{2}}{2(\tau -|\mathbf{x}- 
\mathbf{x}_{0}|\cos \theta )}.  \label{106}
\end{equation}%
We take $\varphi $ as the second curvilinear coordinate and 
\begin{equation}
z=\xi _{3}^{\prime }=r\cos \theta =\frac{(\tau ^{2}-|\mathbf{x}-\mathbf{x}
_{0}|^{2})\cos \theta }{2(\tau -|\mathbf{x}-\mathbf{x}_{0}|\cos \theta )}
,~\theta \in \lbrack 0,\pi ],  \label{107}
\end{equation}%
as the third coordinate.

Express $\cos \theta $ and $\sin \theta $ through $z$ and $\tau $. Using (%
\ref{106}), we obtain 
\begin{eqnarray}
\cos \theta (z,\tau ,\mathbf{x}) &=&\frac{2z\tau }{\tau ^{2}-|\mathbf{x}-%
\mathbf{x}_{0}|^{2}+2z|\mathbf{x}-\mathbf{x}_{0}|},  \notag \\
\sin \theta (z,\tau ,\mathbf{x}) &=&\frac{\big[\big(\tau ^{2}-|\mathbf{x}-%
\mathbf{x}_{0}|^{2}+2z|\mathbf{x}-\mathbf{x}_{0}|\big)^{2}-4z^{2}\tau ^{2}%
\big]^{1/2}}{\tau ^{2}-|\mathbf{x}-\mathbf{x}_{0}|^{2}+2z|\mathbf{x}-\mathbf{%
\ x}_{0}|}.  \label{108}
\end{eqnarray}%
Formula (\ref{106}) allows us to find the equation of the ellipsoid $E(%
\mathbf{x},\mathbf{x}_{0},\tau )$ in coordinates $z,\tau $: 
\begin{equation}
r=r(z,\tau ,\mathbf{x})=\frac{z}{\cos \theta (z,\tau )}=\frac{1}{2\tau }%
\left( \tau ^{2}-|\mathbf{x}-\mathbf{x}_{0}|^{2}+2z|\mathbf{x}-\mathbf{x}%
_{0}|\right) .  \label{109}
\end{equation}%
It follows from (\ref{106}) that 
\begin{equation*}
r\in \lbrack r_{1}(\tau ,\mathbf{x}),r_{2}(\tau ,\mathbf{x})],~r_{1}(\tau ,%
\mathbf{x})=\frac{1}{2}(\tau -|\mathbf{x}-\mathbf{x}_{0}|),~r_{2}(\tau ,%
\mathbf{x})=\frac{1}{2}(\tau +|\mathbf{x}-\mathbf{x}_{0}|).
\end{equation*}%
In addition, (\ref{107}) implies 
\begin{equation*}
z\in \lbrack z_{1}(\tau ,\mathbf{x}),z_{2}(\tau ,\mathbf{x})],~z_{1}(\tau ,%
\mathbf{x})=-\frac{1}{2}(\tau -|\mathbf{x}-\mathbf{x}_{0}|),~z_{2}(\tau .%
\mathbf{x})=\frac{1}{2}(\tau +|\mathbf{x}-\mathbf{x}_{0}|).
\end{equation*}%
Using (\ref{108}), we obtain%
\begin{equation*}
\left. 
\begin{array}{c}
r(z,\tau ,\mathbf{x})\sin \theta =\frac{1}{2\tau }\big[\big(\tau ^{2}-|%
\mathbf{x}-\mathbf{x}_{0}|^{2}+ \\ 
+2z|\mathbf{x}-\mathbf{x}_{0}|\big)^{2}-4z^{2}\tau ^{2}\big]^{1/2}=:f(z,\tau
,\mathbf{x}).%
\end{array}%
\right.
\end{equation*}
Hence, (\ref{105}) becomes 
\begin{equation}
\left. 
\begin{array}{c}
\xi =\mathbf{x}_{0}+\sum_{i=1}^{3}\xi _{i}^{\prime }\mathbf{e}_{i},~\xi
_{1}^{\prime }=f(z,\tau ,\mathbf{x})\cos \varphi ,~ \\ 
\xi _{2}^{\prime }=f(z,\tau ,\mathbf{x})\sin \varphi ,~\xi _{3}^{\prime
}=z,~\varphi \in \lbrack 0,2\pi ).%
\end{array}%
\right.  \label{111}
\end{equation}%
Hence, 
\begin{equation*}
d\xi =d\xi ^{\prime }=f(z,\tau ,\mathbf{x})\frac{\partial f(z,\tau ,\mathbf{x%
})}{\partial \tau }\,d\tau dzd\varphi .
\end{equation*}%
Using 
\begin{equation*}
f(z,\tau ,\mathbf{x})\frac{\partial f(z,\tau ,\mathbf{x})}{\partial \tau }=%
\frac{1}{2}\frac{\partial f^{2}(z,\tau ,\mathbf{x})}{\partial \tau },
\end{equation*}%
we obtain 
\begin{equation*}
\frac{1}{2}\frac{\partial f^{2}(z,\tau ,\mathbf{x})}{\partial \tau }=\frac{1%
}{4\tau ^{3}}\left( \tau ^{2}-|\mathbf{x}-\mathbf{x}_{0}|^{2}+2z|\mathbf{x}-%
\mathbf{x}_{0}|\right) \left( \tau ^{2}+|\mathbf{x}+\mathbf{x}_{0}|^{2}-2z|%
\mathbf{x}-\mathbf{x}_{0}|\right)
\end{equation*}%
Moreover, using $\mathbf{x}-\mathbf{x}_{0}=|\mathbf{x}-\mathbf{x}_{0}|%
\mathbf{e}_{3}$ and (\ref{111}), we obtain 
\begin{equation*}
|\xi -\mathbf{x}|=|\xi -\mathbf{x}_{0}-(\mathbf{x}-\mathbf{x}_{0})|=\left(
r^{2}(\tau ,z)-2z|\mathbf{x}-\mathbf{x}_{0}|+|\mathbf{x}-\mathbf{x}%
_{0}|^{2}\right) ^{1/2}.
\end{equation*}%
Substitute in this formula $r(\tau ,z,\mathbf{x})$ from (\ref{109}). After
some calculations, we obtain 
\begin{equation*}
|\xi -\mathbf{x}|=\frac{\tau ^{2}+|\mathbf{x}-\mathbf{x}_{0}|^{2}-2z|\mathbf{%
\ x}-\mathbf{x}_{0}|}{2\tau }.
\end{equation*}%
Hence, 
\begin{eqnarray}
\frac{d\xi }{|\xi -\mathbf{x}|} &=&\frac{1}{2\tau ^{2}}\left( \tau ^{2}-|%
\mathbf{x}-\mathbf{x}_{0}|^{2}+2z|\mathbf{x}-\mathbf{x}_{0}|\right) \,d\tau
dzd\varphi  \label{112} \\
&=&\frac{r(z,\tau ,\mathbf{x})}{\tau }\,d\tau dzd\varphi =\frac{|\xi -%
\mathbf{x}_{0}|}{\tau }\,d\tau dzd\varphi .  \notag
\end{eqnarray}%
Taking into account (\ref{112}), rewrite equation (\ref{104}) as%
\begin{equation}
u(\mathbf{x},\mathbf{x}_{0},t)=\frac{1}{4\pi |\mathbf{x}-\mathbf{x}_{0}|}+
\label{113}
\end{equation}%
\begin{equation*}
+\frac{1}{4\pi }\int\limits_{|\mathbf{x}-\mathbf{x}_{0}|}^{t}\int%
\limits_{z_{1}(\tau ,\mathbf{x})}^{z_{2}(\tau ,\mathbf{x})}\int%
\limits_{0}^{2\pi }a(\xi ,t-|\xi -\mathbf{x}|)u(\xi ,\mathbf{x}_{0},t-|\xi -%
\mathbf{x}|)\frac{|\xi -\mathbf{x}_{0}|}{\tau }\,d\varphi dzd\tau ,(\mathbf{x%
},t)\in K_{\mathbf{x}_{0},T}^{+},
\end{equation*}%
where $\xi $ is given in (\ref{111}).

Rewrite (\ref{113}) as the equation for function $\overline{u}(\mathbf{x},%
\mathbf{x}_{0},t)$,%
\begin{equation*}
\overline{u}(\mathbf{x},\mathbf{x}_{0},t)=\frac{1}{16\pi ^{2}}\int\limits_{|%
\mathbf{x}-\mathbf{x}_{0}|}^{t}\int\limits_{z_{1}(\tau ,\mathbf{x}%
)}^{z_{2}(\tau ,\mathbf{x})}\int\limits_{0}^{2\pi }a(\xi ,t-|\xi -\mathbf{x}%
|)\frac{1}{\tau }\,d\varphi dzd\tau +
\end{equation*}%
\begin{equation}
+\frac{1}{4\pi }\int\limits_{|\mathbf{x}-\mathbf{x}_{0}|}^{t}\int%
\limits_{z_{1}(\tau ,\mathbf{x})}^{z_{2}(\tau ,\mathbf{x})}\int%
\limits_{0}^{2\pi }a(\xi ,t-|\xi -\mathbf{x}|)\overline{u}(\xi ,\mathbf{x}%
_{0},t-|\xi -\mathbf{x}|)\frac{|\xi -\mathbf{x}_{0}|}{\tau }\,d\varphi
dzd\tau ,  \label{114}
\end{equation}%
\begin{equation*}
(\mathbf{x},t)\in K_{\mathbf{x}_{0},T}^{+}.
\end{equation*}%
It is convenient to transform this equation in a different one. Replace the
variable $z$ by $z^{\prime }$as 
\begin{equation*}
z=\frac{1}{2}(|\mathbf{x}-\mathbf{x}_{0}|+\tau z^{\prime }).
\end{equation*}%
Hence, 
\begin{equation*}
\frac{dz}{\tau }=\frac{1}{2}dz^{\prime }.
\end{equation*}%
Hence, 
\begin{equation*}
f(z,\tau ,\mathbf{x})|_{z=(|\mathbf{x}-\mathbf{x}_{0}|+z^{\prime }\tau
)/2}=:f^{\ast }(z^{\prime },\tau ,\mathbf{x})=\frac{1}{2}\tau (\tau ^{2}-|%
\mathbf{x}-\mathbf{x}_{0}|^{2})(1-(z^{^{\prime }})^{2},
\end{equation*}%
We have%
\begin{equation}
\left. 
\begin{array}{c}
\xi =\xi ^{\ast }(z^{\prime },\tau ,\mathbf{x})=\mathbf{x}%
_{0}+\sum_{i=1}^{3}\xi _{i}^{\prime }\mathbf{e}_{i},~ \\ 
\xi _{1}^{\prime }=f^{\ast }(z^{\prime },\tau ,\mathbf{x})\cos \varphi ,~\xi
_{2}^{\prime }=f^{\ast }(z^{\prime },\tau ,\mathbf{x})\sin \varphi , \\ 
\xi _{3}^{\prime }=z=\left( 1/2\right) (|\mathbf{x}-\mathbf{x}%
_{0}|+z^{\prime }\tau ),~\varphi \in \lbrack 0,2\pi ).%
\end{array}%
\right.  \label{115}
\end{equation}

Then integral equation (\ref{114}) becomes%
\begin{equation}
\overline{u}(\mathbf{x},\mathbf{x}_{0},t)=\frac{1}{32\pi ^{2}}\int\limits_{| 
\mathbf{x}-\mathbf{x}_{0}|}^{t}\int\limits_{-1}^{1}\int\limits_{0}^{2\pi
}[a(\xi ,t-|\xi -\mathbf{x}|)]_{\xi =\xi ^{\ast }(z^{\prime },\tau ,\mathbf{%
x })}\,d\varphi dz^{\prime }d\tau +  \label{116}
\end{equation}%
\begin{equation*}
+\frac{1}{8\pi }\int\limits_{|\mathbf{x}-\mathbf{x}_{0}|}^{t}\int
\limits_{-1}^{1}\int\limits_{0}^{2\pi }[a(\xi ,t-|\xi -\mathbf{x}|)\overline{
u}(\xi ,\mathbf{x}_{0},t-|\xi -\mathbf{x}|)|\xi -\mathbf{x}_{0}|]_{\xi =\xi
^{\ast }(z^{\prime },\tau ,\mathbf{x})}\,d\varphi dz^{\prime }d\tau ,
\end{equation*}%
\begin{equation*}
(\mathbf{x},t)\in K_{\mathbf{x}_{0},T}^{+}.
\end{equation*}

We solve integral equation (\ref{116}) by the method of successive
approximations as:%
\begin{equation*}
\overline{u}_{1}(\mathbf{x},\mathbf{x}_{0},t)=\frac{1}{16\pi ^{2}}
\int\limits_{|\mathbf{x}-\mathbf{x}_{0}|}^{t}\int\limits_{-1}^{1}\int
\limits_{0}^{2\pi }[a(\xi ,t-|\xi -\mathbf{x}|)]_{\xi =\xi ^{\ast
}(z^{\prime },\tau ,\mathbf{x})}\,d\varphi dz^{\prime }d\tau ,
\end{equation*}%
\begin{equation}
\overline{u}_{n}(\mathbf{x},\mathbf{x}_{0},t)=\frac{1}{8\pi }\int\limits_{| 
\mathbf{x}-\mathbf{x}_{0}|}^{t}\int\limits_{-1}^{1}\int\limits_{0}^{2\pi
}[a(\xi ,t-|\xi -\mathbf{x}|)\times  \label{117}
\end{equation}%
\begin{equation*}
\times \overline{u}_{n-1}(\xi ,\mathbf{x}_{0},t-|\xi -\mathbf{x}|)|\xi - 
\mathbf{x}_{0}|]|_{\xi =\xi ^{\ast }(z^{\prime },\tau ,\mathbf{x}
)}\,d\varphi dz^{\prime }d\tau ,
\end{equation*}%
\begin{equation*}
n=2,3,\ldots ,~(\mathbf{x},t)\in K_{\mathbf{x}_{0},T}^{+}.
\end{equation*}

Let 
\begin{equation*}
\max_{(\mathbf{x},t)\in K_{\mathbf{x}_{0},T}^{+}}a(\mathbf{x},t)=A_{0}.
\end{equation*}%
Then 
\begin{equation*}
0\leq \overline{u}_{1}(\mathbf{x},\mathbf{x}_{0},t)\leq \frac{1}{32\pi ^{2}}
\int\limits_{|\mathbf{x}-\mathbf{x}_{0}|}^{t}\int\limits_{-1}^{1}\int
\limits_{0}^{2\pi }A_{0}\,d\varphi dz^{\prime }d\tau \leq \frac{A_{0}}{8\pi }
(t-|(\mathbf{x}-\mathbf{x}_{0}|),~(\mathbf{x},t)\in K_{\mathbf{x}_{0},T}^{+}.
\end{equation*}%
Next, we find 
\begin{eqnarray*}
0 &\leq &\overline{u}_{2}(\mathbf{x},\mathbf{x}_{0},t)\leq \frac{A_{0}}{8\pi 
}\int\limits_{|\mathbf{x}-\mathbf{x}_{0}|}^{t}\int\limits_{-1}^{1}\int
\limits_{0}^{2\pi }\frac{A_{0}(t-|\xi ^{\ast }-\mathbf{x}|-|\xi ^{\ast }- 
\mathbf{x}_{0}|)}{8\pi \tau }|\xi ^{\ast }-\mathbf{x}_{0}|\,d\varphi
dz^{\prime }d\tau \leq \\
&\leq &\frac{A_{0}^{2}T}{16\pi }\int\limits_{|\mathbf{x}-\mathbf{x}
_{0}|}^{t}(t-\tau )\,d\tau =\frac{A_{0}^{2}}{16\pi }\frac{(t-|\mathbf{x}- 
\mathbf{x}_{0}|)^{2}}{2!},~(\mathbf{x},t)\in K_{\mathbf{x}_{0},T}^{+}.
\end{eqnarray*}%
We have used above the inequality $|\xi ^{\ast }-\mathbf{x}_{0}|\leq T$.
Using the mathematical induction method, we obtain 
\begin{equation}
0\leq \overline{u}_{n}(\mathbf{x},\mathbf{x}_{0},t)\leq \frac{
A_{0}^{n}T^{n-1}}{4\cdot 2^{n}\pi }\frac{(t-|\mathbf{x}-\mathbf{x}_{0}|)^{n} 
}{n!},~(\mathbf{x},t)\in K_{\mathbf{x}_{0},T}^{+},\text{ }n=1,2,\ldots
\label{1}
\end{equation}%
Note that all functions $\overline{u}_{n}(\mathbf{x},\mathbf{x}_{0},t)$, $%
n=1,2,\ldots ,$ are continuous in $\overline{K_{\mathbf{x}_{0},T}^{+}}$.
Since 
\begin{eqnarray}
0 &\leq &\overline{u}(\mathbf{x},\mathbf{x}_{0},t)=\sum\limits_{n=1}^{\infty
}\overline{u}_{n}(\mathbf{x},\mathbf{x}_{0},t)\leq \sum\limits_{n=1}^{\infty
}\frac{A_{0}^{n}T^{n-1}}{4\cdot 2^{n}\pi }\frac{(t-|\mathbf{x}-\mathbf{x}
_{0}|)^{n}}{n!}=  \notag \\
&=&\frac{1}{4\pi T}\left[ \exp \left( \frac{A_{0}T(t-|\mathbf{x}-\mathbf{x}
_{0}|)}{2}\right) -1\right]  \label{2} \\
&\leq &\frac{1}{4\pi T}\left[ \exp \left( \frac{A_{0}T^{2}}{2}\right) -1 %
\right] ,~(\mathbf{x},t)\in K_{\mathbf{x}_{0},T}^{+},  \notag
\end{eqnarray}%
and the latter series converges uniformly in $K_{\mathbf{x}_{0},T}^{+}$,
then (\ref{1}) and (\ref{2}) imply that the function $\overline{u}(\mathbf{x}%
,\mathbf{x}_{0},t)$ is nonnegative and continuous in $K_{\mathbf{x}
_{0},T}^{+}$. Hence, inequality (\ref{102}) is true and $\overline{u}\in C(%
\overline{K_{\mathbf{x}_{0},T}^{+}})$. In particular, (\ref{201}) implies
that $\overline{u}\in C\left( \overline{Q\left( R,T^{-},T\right) }\right) .$

We now prove that the function $\overline{u}(\mathbf{x},\mathbf{x}_{0},t)\in
C^{2}\left( \overline{K_{,\mathbf{x}_{0},T}^{+}}\right) \cap (\overline{
\Omega }\times \lbrack 0,T])$, which implies that $\overline{u}\in
C^{2}\left( \overline{Q\left( R,T^{-},T\right) }\right) ,$ see (\ref{201}).
We start calculations with the derivative with respect to $t$.
Differentiating equation (\ref{116}), we obtain 
\begin{equation*}
\overline{u}_{t}(\mathbf{x},\mathbf{x}_{0},t)=\frac{1}{32\pi }
\int\limits_{-1}^{1}\int\limits_{0}^{2\pi }[a(\xi ,t-|\xi -\mathbf{x}|)]\mid
_{\xi =\xi ^{\ast }(z^{\prime },t,\mathbf{x})}\,d\varphi dz^{\prime }+
\end{equation*}%
\begin{equation*}
+\frac{1}{32\pi }\int\limits_{|\mathbf{x}-\mathbf{x}_{0}|}^{t}\int
\limits_{-1}^{1}\int\limits_{0}^{2\pi }[a_{t}(\xi ,t-|\xi -\mathbf{x}|)]\mid
_{\xi =\xi ^{\ast }(z^{\prime },\tau ,\mathbf{x})}\,d\varphi dz^{\prime
}d\tau +
\end{equation*}%
\begin{equation}
+\frac{1}{8\pi }\int\limits_{|\mathbf{x}-\mathbf{x}_{0}|}^{t}\int
\limits_{-1}^{1}\int\limits_{0}^{2\pi }[a_{t}(\xi ,t-|\xi -\mathbf{x}|) 
\overline{u}(\xi ,t-|\xi -\mathbf{x}|)|\xi -\mathbf{x}_{0}|]\mid _{\xi =\xi
^{\ast }(z^{\prime },\tau ,\mathbf{x})}\,d\varphi dz^{\prime }d\tau +
\label{118}
\end{equation}%
\begin{equation*}
+\frac{1}{8\pi }\int\limits_{|\mathbf{x}-\mathbf{x}_{0}|}^{t}\int
\limits_{_{1}}^{1}\int\limits_{0}^{2\pi }a(\xi ,t-|\xi -\mathbf{x}|) 
\overline{u}_{t}(\xi ,\mathbf{x}_{0},t-|\xi -\mathbf{x}|)|\xi -\mathbf{x}
_{0}|]\mid _{\xi =\xi ^{\ast }(z^{\prime },\tau ,\mathbf{x})}\,d\varphi
dz^{\prime }d\tau ,~
\end{equation*}%
\begin{equation*}
(\mathbf{x},t)\in K_{\mathbf{x}_{0},T}^{+}.
\end{equation*}%
The first four integrals in (\ref{118}) are continuous known functions in $%
\overline{K_{\mathbf{x}_{0},T}^{+}}$. Therefore (\ref{118}) present the
equation quite similar to equation (\ref{116}). Thus. applying the method of
successive approximations, we can prove that $\overline{u}_{t}\in C(%
\overline{K_{\mathbf{x}_{0},T}^{+}})$.

We remind that $K_{\mathbf{x}_{0},T}^{+}$ is defined in (\ref{200}). We now
show that partial derivatives with respect to space variables are continuous
in domain $Q\left( R,T^{-},T\right) $. However, they are not continuous in $%
K_{\mathbf{x}_{0},T}^{+}\diagdown (\overline{\Omega }\times \lbrack 0,T])$
since the derivatives of $|\mathbf{x}-\mathbf{x}_{0}|$ are not continuous in
a vicinity of the point $\mathbf{x}_{0}$. Using 
\begin{equation*}
|\xi ^{\ast }(z^{\prime },\tau ,\mathbf{x})-\mathbf{x}|=\tau -|\xi ^{\ast
}(z^{\prime },\tau ,\mathbf{x})-\mathbf{x}_{0}|,
\end{equation*}%
we rewrite equation (\ref{116}) as%
\begin{equation*}
\overline{u}(\mathbf{x},\mathbf{x}_{0},t)=\frac{1}{32\pi ^{2}}\int\limits_{| 
\mathbf{x}-\mathbf{x}_{0}|}^{t}\int\limits_{-1}^{1}\int\limits_{0}^{2\pi
}[a(\xi ,t-\tau +|\xi -\mathbf{x}_{0}|)]\mid _{\xi =\xi ^{\ast }(z^{\prime
},\tau ,\mathbf{x})}\,d\varphi dz^{\prime }d\tau +
\end{equation*}%
\begin{equation}
+\frac{1}{8\pi }\int\limits_{|\mathbf{x}-\mathbf{x}_{0}|}^{t}\int
\limits_{-1}^{1}\int\limits_{0}^{2\pi }[a(\xi ,t-\tau +|\xi -\mathbf{x}
_{0}|)|\xi -\mathbf{x}_{0}|\times  \label{119}
\end{equation}%
\begin{equation*}
\times \overline{u}(\xi ,\mathbf{x}_{0},t-\tau +|\xi -\mathbf{x}_{0}|)]_{\xi
=\xi ^{\ast }(z^{\prime },\tau ,\mathbf{x})}\,d\varphi dz^{\prime }d\tau ,~( 
\mathbf{x},t)\in K_{\mathbf{x}_{0},T}^{+}.
\end{equation*}

We now differentiate equation (\ref{119}) with respect to $x_{k}$, $k=1,2,3$
setting $x_{1}=x,x_{2}=y,x_{3}=z$. We obtain 
\begin{eqnarray}
\partial _{x_{k}}\overline{u}(\mathbf{x},\mathbf{x}_{0},t) &=&-\frac{
\partial _{x_{k}}|\mathbf{x}-\mathbf{x}_{0}|}{32\pi }\int\limits_{-1}^{1}
\int\limits_{0}^{2\pi }[a(\xi ,t-|\mathbf{x}-\mathbf{x}_{0}|+|\xi -\mathbf{x}
_{0}|)]\mid _{\xi =\xi ^{\ast \ast }(z^{\prime },\mathbf{x})}\,d\varphi
dz^{\prime }-  \notag \\
&&-\frac{\partial _{x_{k}}|\mathbf{x}-\mathbf{x}_{0}|}{8\pi }
\int\limits_{-1}^{1}\int\limits_{0}^{2\pi }[a(\xi ,t-|\mathbf{x}-\mathbf{x}
_{0}|+|\xi -\mathbf{x}_{0}|)|\xi -\mathbf{x}_{0}|\times  \label{121} \\
&&\times \overline{u}(\xi ,\mathbf{x}_{0},t-|\mathbf{x}-\mathbf{x}_{0}|+|\xi
-\mathbf{x}_{0}|)]\mid _{\xi =\xi ^{\ast \ast }(z^{\prime },\mathbf{x}
)}\,d\varphi dz^{\prime }+  \notag \\
&&+\frac{1}{32\pi }\int\limits_{|\mathbf{x}-\mathbf{x}_{0}|}^{t}\int
\limits_{-1}^{1}\int\limits_{0}^{2\pi }\{\nabla _{\xi }[a(\xi ,t-\tau +|\xi
- \mathbf{x}_{0}|)]_{\xi =\xi ^{\ast }(z^{\prime },\tau ,\mathbf{x})}\cdot
\partial _{x_{k}}\xi ^{\ast }(z^{\prime },\tau ,\mathbf{x})+  \notag \\
&&+[a_{t}(\xi ,t-\tau +|\xi -\mathbf{x}_{0}|)]_{\xi =\xi ^{\ast }(z^{\prime
},\tau ,\mathbf{x})}\partial _{x_{k}}|\xi ^{\ast }(z^{\prime },\tau ,\mathbf{%
\ x})-\mathbf{x}_{0}|\}\,d\varphi dz^{\prime }d\tau +  \notag \\
&&+\frac{1}{8\pi }\int\limits_{|\mathbf{x}-\mathbf{x}_{0}|}^{t}\int
\limits_{-1}^{1}\int\limits_{0}^{2\pi }\{\nabla _{\xi }[a(\xi ,t-\tau +|\xi
- \mathbf{x}_{0}|)|\xi -\mathbf{x}_{0}|]_{\xi =\xi ^{\ast }(z^{\prime },\tau
, \mathbf{x})}\cdot \partial _{x_{k}}\xi ^{\ast }(z^{\prime },\tau ,\mathbf{x%
})  \notag \\
&&+[a_{t}(\xi ,t-\tau +|\xi -\mathbf{x}_{0}|)|\xi -\mathbf{x}_{0}|]_{\xi
=\xi ^{\ast }(z^{\prime },\tau ,\mathbf{x})}\partial _{x_{k}}|\xi ^{\ast
}(z^{\prime },\tau ,\mathbf{x})-\mathbf{x}_{0}|\}  \notag \\
&&\times \lbrack \overline{u}(\xi ,\mathbf{x}_{0},t-\tau +|\xi -\mathbf{x}
_{0}|)]_{\xi =\xi ^{\ast }(z^{\prime },\tau ,\mathbf{x})}\,d\varphi
dz^{\prime }d\tau  \notag \\
&&+\frac{1}{8\pi }\int\limits_{|\mathbf{x}-\mathbf{x}_{0}|}^{t}\int
\limits_{_{1}}^{1}\int\limits_{0}^{2\pi }[a(\xi ,t-\tau +|\xi -\mathbf{x}
_{0}|)|\xi -\mathbf{x}_{0}|]_{\xi =\xi ^{\ast }(z^{\prime },\tau ,\mathbf{x}
)}\times  \notag \\
&&\times \{\nabla _{\xi }[\overline{u}(\xi ,\mathbf{x}_{0},t-\tau +|\xi - 
\mathbf{x}_{0}|)]_{\xi =\xi ^{\ast }(z^{\prime },\tau ,\mathbf{x})}\cdot
\partial _{x_{k}}\xi ^{\ast }(z^{\prime },\tau ,\mathbf{x})+  \notag \\
&&+[\overline{u}_{t}(\xi ,\mathbf{x}_{0},t-\tau +|\xi -\mathbf{x}
_{0}|)]_{\xi =\xi ^{\ast }(z^{\prime },\tau ,\mathbf{x})}\partial
_{x_{k}}|\xi ^{\ast }(z^{\prime },\tau ,\mathbf{x})-\mathbf{x}
_{0}|\}\,d\varphi dz^{\prime }d\tau ,  \notag \\
k &=&1,2,3,~(\mathbf{x},t)\in K_{\mathbf{x}_{0},T}^{+}\cap (\overline{\Omega 
}\times \lbrack 0,T]).  \notag
\end{eqnarray}%
where 
\begin{equation*}
\xi ^{\ast \ast }(z^{\prime },\mathbf{x})=\xi ^{\ast }(z^{\prime },|\mathbf{%
x }-\mathbf{x}_{0}|,\mathbf{x})=\mathbf{x}_{0}+\frac{1}{2}(\mathbf{x}-%
\mathbf{x }_{0})(1+z^{\prime }).
\end{equation*}

Three equations (\ref{121}) for $k=1,2,3$ form the system of three integral
equations of the Volterra type. This system similar to the scalar equation (%
\ref{116}), and it is closed in the domain $(\mathbf{x},t)\in K_{\mathbf{x}%
_{0},T}^{+}\cap (\overline{\Omega }\times \lbrack 0,T])$. Kernels of this
equations are continuous functions since $|\mathbf{x}-\mathbf{x}_{0}|\geq 
\text{dist}(\mathbf{x}_{0},\Omega )>0$, if $(\mathbf{x},t)\in K_{\mathbf{x}%
_{0},T}^{+}\cap (\overline{\Omega }\times \lbrack 0,T]).$ Hence, one can
apply the method of successive approximation to prove existence and
uniqueness of a continuous solution of equations (\ref{121}).

Acting similarly via calculating the second order derivatives $\overline{u}%
_{tt}$, $\partial _{x_{k}}\overline{u}_{t}$, $\partial _{x_{k}x_{j}}%
\overline{u}$, one can obtain another system of integral equations of the
Volterra type and prove that all these derivatives are continuous in $(%
\mathbf{x},t)\in K_{\mathbf{x}_{0},T}^{+}\cap (\overline{\Omega }\times
\lbrack 0,T]).$ $\square $\hfill $\Box $




\end{document}